\documentclass[11pt]{amsart}
\usepackage{amssymb,amsthm,amsfonts,amsmath}
\usepackage[pdftex]{graphicx}
\usepackage[mathscr]{euscript}
\usepackage[utf8]{inputenc}
\usepackage[english]{babel}
\newtheorem{theorem}{Theorem}
\newtheorem{lemma}{Lemma}

\newtheorem{remark}{Remark}

\theoremstyle{definition}

\begin{document}
	
	\title[Decomposition of some classes of regular graphs into ...]{Decomposition of some classes of regular graphs into cycles and paths of length eight}
	
	\author{Cecily Sahai. C}
	\address{Department of Mathematics, Sri Sivasubramaniya Nadar College of Engineering, Chennai, India, postcode-603110}
	\email{cecilysahai@gmail.com}
	
	\author{Sampath Kumar. S}
	\address{Department of Mathematics, Sri Sivasubramaniya Nadar College of Engineering, Chennai, India, postcode-603110}
	\email{sampathkumars@ssn.edu.in}
	
	\author{Arputha Jose. T}
	\address{Department of Mathematics, Sri Sivasubramaniya Nadar College of Engineering, Chennai, India, postcode-603110}
	\email{arputhajose792@gmail.com}
	
	%

	\thanks{}
	
	\begin{abstract}
		In this paper, we examine the necessary and sufficient conditions for the existence of a $(8; p, q)$-decomposition of $\lambda$-fold tensor product and wreath product of complete graphs.
	\end{abstract}
\keywords{Cycle decomposition, Path decomposition, Tensor product, Wreath product}

	\maketitle
	
	\section{Introduction}
	\indent We consider only finite undirected graphs. For a graph $G$, $G(\lambda)$ represent the graph obtained by replacing every edge of $G$ by $\lambda$ parallel edge and $\lambda G$ represent  $\lambda$ edge-disjoint copies of $G$. Let $C_{k}$ (resp. $P_{k}$) denote the cycle (resp. path) of length $k$. The complete graph on $n$ vertices is denoted by $K_{n}$. Let $\langle X \rangle$ denotes subgraph induced by $X$. Let $\langle X,Y \rangle$ denote the set of all edges with one end in $X$ and other end in $Y$. If $H_{1}, H_{2},\ldots,H_{k}$ are edge disjoint subgraphs of $G$ such that $E\left( G \right) = E\left( H_{1}\right) \cup E\left( H_{2}\right) \cup \ldots \cup E\left( H_{k}\right)$, then we say that $H_{1}, H_{2},\ldots,H_{k}$ decompose $G$ which is denoted by $G=H_{1} \oplus H_{2}\oplus\ldots\oplus H_{k}$. When each $H_{i} \cong H$, then we say that $G$ has a $H$-decomposition which is denoted by $H|G$. When $G$ can be decomposed into $p$ copies of $ P_{k}$ and $q$ copies of $ C_{k}$ such that $k(p+q)=E(G)$, then we say $G$ admits a $(k; p, q)$-decomposition. A $k$-factor of a graph $G$ is a $k$-regular spanning subgraph, which is denoted by $k$F and $G-F_{k}$ denote a graph after the removal of a $k$-factor. For two graphs $G$ and $H,$ their wreath product $G\circ H$ has vertex set $V(G)\times V(H) $ in which two vertices $\left( g_{1},h_{1}\right)$ and $\left( g_{2},h_{2}\right)$ are adjacent whenever $g_{1}g_{2}$ is an edge in $G$ or $g_{1}=g_{2}$ and $h_{1}h_{2}$ is an edge in $H$. Similarly, the tensor product $G \times H$ has vertex set $V (G) \times V (H)$ in which $(g_{1}, h_{1})(g_{2}, h_{2})$ is an edge whenever $g_{1}g_{2}$ is an edge in $G$ and $h_{1}h_{2}$ is an edge in $H$. For the standard graph-theoretic terminologies, the reader is referred to \cite{TBGT2012} and \cite{GTGTM2008}.\\
	
	\indent In \cite{DCGPC2010}, Shyu obtained neccessary and sufficient conditions for decomposition of complete graphs into paths and cycles of length four. The necessary conditions for the existence of a decomposition of complete equipartite graph (3,4 and 5 parts) into paths and cycles are proved to be sufficient in \cite{PCDCEGFP2009} and \cite{PCDCEG35P2010}. In \cite{DCBGPC2014}, Jeevadoss et al. obtained neccessary and sufficient conditions for existence of $(k; p, q)$-decomposition of $K_{m,n}$ and $K_{n}$, when $m \geq \frac{k}{2}, n \geq \lceil \frac{(k+1)}{2} \rceil$ for $k \equiv 0($mod $4)$ and when $m,n \geq 2k$ for $k \equiv 2$(mod $4$). The same authors obtained $(4; p, q)$-decomposition of tensor product, wreath product and cartesian product of complete graphs in \cite{DPGPCLF2016}. In \cite{DCTGCP32020}, decomposition of complete tripartite graphs into cycles and paths of length three is established. In \cite{DHGPC2022}, neccessary and sufficient conditions for decomposition of hypercube graphs into paths and cycles of length four is obtained. In \cite{DCTGCP2023}, decomposition of complete tripartite graphs into cycles and paths of prime length is established. In \cite{DCEGPC2P2023}, Pauline Ezhilarasi et al. obtained neccessary and sufficient conditions for decomposition of complete equipartite graphs into paths and cycles of size $2p$. In \cite{DLGPC2023}, Ganesamurthy gave neccessary and sufficient conditions for existence of $(k; p, q)$-decomposition of line graph of complete graph when $n \equiv 0,1,2($mod $k),\ k \geq 7$ and $k$ is prime.\\
	In this paper, we study the $(8; p, q)$-decomposition of $(K_{m} \times K_{n})(\lambda)$ and $(K_{m}\circ \overline{K}_{n})(\lambda)$.  \\
	We quote some known results for our future reference.
	
	\begin{theorem} [\cite{DCMSPNHD1983}]
		A necessary and sufficient condition for the existence of decomposition of a complete multigraph $K_{n}(\lambda)$ into edge disjoint simple paths of length $ m $ is $\lambda n(n-1) \equiv 0($mod $2m)$ and $n \geq m+1$
		\label{THM1C8P8}	
	\end{theorem}
	
	\begin{theorem} [\cite{DT2009}]
		Let $  m \geq 3 $ be an odd integer.(1) If $  m \equiv 1$ or $3($mod $ 6) $, then $ K_{3} | K_{m} $. (2) If $ m \equiv 5($mod $ 6) $, then $ K_{m} $ can be decomposed into $ \left(\frac{m(m-1)-20}{6}\right) \ K_{3}$ and a $K_{5}$.
		\label{THM2C8P8}	
	\end{theorem}
	
	\begin{theorem} [\cite{NDP1985}]
		Let $m \geq n$ and either $\lambda$ is even or $m$ and $n$ are even. Then $K_{m,n}(\lambda)$ is decomposed into $P_{k}$ if and only if $\lambda mn \equiv 0($mod $k-1)$, $m \geq \lceil \frac{k}{2} \rceil$ and $n \geq \lceil \frac{k-1}{2} \rceil$.
		\label{THM3C8P8}
	\end{theorem}
	
	\begin{theorem} [\cite{CDCM2011}]
		Let $ \lambda, n $ and $ m $ be integers with $ n, m \geq 3 $ and $ \lambda \geq 1 $. There exists adecomposition of $K_{n}(\lambda)$ into $ m $-cycles if and only if (1) $m \leq n$ (2) $\lambda (n-1)$ is even (3) $m$ divides $\lambda {n \choose 2}$
		\label{THM5C8P8}	
	\end{theorem}
	
	\begin{theorem} [\cite{DCBMCS2015}]
		For positive integers $ \lambda , k, m $ and $ n $ with $ \lambda m \equiv \lambda n \equiv k \equiv 0 ($mod $2)$ and $\min\{m, n\} \geq k/2 \geq 2 $, the multigraph $K_{m,n}(\lambda)$ is $ C_{k} $-decomposable if one of the following conditions holds: (1) $\lambda$ is odd and $ k $ divides $ mn $, (2) $\lambda$ is even and $ k $ divides $ 2mn $, (3) $\lambda$ is even and $\lambda n$ or $\lambda m$ is divisible by $k$.
		\label{THM6C8P8}	
	\end{theorem}
	
	\section{Prerequisites}
	In this section, we will prove the $(8; p, q)$-decomposition of some small graphs.\\
	\begin{remark}
		If there are $t$ odd degree vertices in $G$, then any $(8; p, q)$-decomposition of $G$ must have atleast $\frac{t}{2}$ paths. 
	\end{remark}
	
	\begin{lemma}
		If $G$ has a $(8; p, q)$-decomposition, then $G \times H$ has a $(8; p, q)$-decomposition.
		\label{LEM0C8P8}
	\end{lemma}
	\noindent{\bf Proof.}
	The graph $G \times H \cong G \times (K_{2} \oplus K_{2} \oplus \ldots K_{2}) \cong (G \times K_{2}) \oplus (G \times K_{2}) \oplus \ldots \oplus (G \times K_{2}) $. Since $G$ has a $(8; p, q)$-decomposition, it is straight forward to check that $G \times K_{2}$ has a $(8; p, q)$-decomposition. Hence $G \times H$ has a $(8; p, q)$-decomposition.\hfill$\Box$	 
	
	\begin{lemma}
		There exists a $(8; p, q)$-decomposition of $K_{4} \times K_{4}$.
		\label{LEM1C8P8}
	\end{lemma}
	\noindent{\bf Proof.}
	Let the vertices of the graph $K_{4} \times K_{4}$ be partioned into $A=\{a_{i}|0 \leq i \leq 3\},\ B=\{b_{i}|0 \leq i \leq 3\},\ C=\{c_{i}|0 \leq i \leq 3\}$ and $D=\{d_{i}|0 \leq i \leq 3\}$. We exhibit a $(8; p, q)$-decomposition of $K_{4} \times K_{4}$ as follows. Let $P=b_{2}c_{1}d_{2}c_{0}a_{3}d_{0}b_{1}c_{2}d_{1}$ and $P^{'}=a_{0}b_{1}d_{2}a_{3}d_{1}b_{3}c_{1}a_{2}c_{0}$. Consider the permutation $\rho=(a_{0}a_{1}a_{2}a_{3})(b_{0}b_{1}b_{2}b_{3})(c_{0}c_{1}c_{2}c_{3})(d_{0}d_{1}d_{2}d_{3})$. Then $\{\rho^{i}(P),\ \rho^{i}(P^{'})| 0 \leq i \leq 3\}$ provides $8$ copies of $ P_{8}$. The remaining edge form an $8$-cycle, say $C=(b_{0}a_{2}b_{1}a_{3}b_{2}a_{0}b_{3}a_{1}b_{0})$. Now, consider $P^{'}\cup C$ and decompose them into $2$ copies of $ P_{8}$ as follows: $a_{0}b_{1}d_{2}a_{3}d_{1}b_{3}c_{1}a_{2}b_{0}$ and $ c_{0}a_{2}b_{1}a_{3}b_{2}a_{0}b_{3}a_{1}b_{0}$. This provides the required $(8; p, q)$-decomposition of $K_{4} \times K_{4}$. \hfill$\Box$
	
	\begin{lemma}
		The graph $P_{2} \circ \overline{K}_{4}$ has a $(8; p, q)$-decomposition.
		\label{LEM2C8P8}
	\end{lemma}
	\noindent{\bf Proof.}
	Let the vertices of the graph $P_{2} \circ \overline{K}_{4}$ be partitioned into $A=\{a_{i}|0 \leq i \leq 3\},\ B=\{b_{i}|0 \leq i \leq 3\}$ and $ C=\{c_{i}|0 \leq i \leq 3\}$. The decomposition of $P_{2} \circ \overline{K}_{4}$ into $8$-cycles given by: $C^{1}= (a_{0} b_{3} a_{3} b_{2} a_{2} b_{1} a_{1} b_{0}),\ C^{2}=(a_{0} b_{1} a_{3} b_{0} a_{2} b_{3} a_{1} b_{2}),\ C^{3}= (b_{0} c_{3} b_{3} c_{2} b_{2} c_{1} b_{1} c_{0})$ and $C^{4}=(b_{0} c_{1} b_{3} c_{0} b_{2} c_{3} b_{1} c_{2})$. The required $(8; p, q)$-decomposition of $P_{2} \circ \overline{K}_{4}$ is as follows:\\
	$D_{1}$: Consider $C^{1}\cup C^{3}$ and decompose them into $2$ copies of $ P_{8}$ as follows: $c_{3} b_{3} a_{3} b_{2} a_{2} b_{1} a_{1} b_{0} a_{0}$ and $a_{0} b_{3} c_{2} b_{2} c_{1} b_{1} c_{0} b_{0} c_{3}$.\\
	$D_{2}$: Decomposition of $C^{2}\cup C^{4}$ into $2$ copies of $ P_{8}$ given by: $c_{3} b_{1} a_{3} b_{0} a_{2} b_{3} a_{1} b_{2} a_{0}$ and $a_{0} b_{1} c_{2} b_{0} c_{1} b_{3} c_{0} b_{2} c_{3}$.\\
	$D_{3}$: By decomposing $C^{1}\cup C^{2}\cup C^{4}$ we obtain $3$ copies of $ P_{8}$ namely, $c_{1} b_{0} a_{3} b_{1} a_{0}\\ b_{2} a_{1} b_{3} a_{2},\ a_{2} b_{0} c_{2} b_{1} c_{3} b_{2} c_{0} b_{3} a_{3}$ and $a_{3} b_{2} a_{2} b_{1} a_{1} b_{0} a_{0} b_{3} c_{1}$.\\
	Here $(p,q) \in \{(0,4),(2,2),(3,1),(4,0)\}$. For $1 \leq i \leq 3$, combinations of $D_{i}$ provide the required $(8; p, q)$-decomposition. \hfill$\Box$
	
	\begin{lemma}
		There exists a $(8; p, q)$-decomposition of $K_{8,8}$.
		\label{LEM3C8P8}
	\end{lemma}
	\noindent{\bf Proof.}
	Partition the vertices of the graph  $K_{8,8}$ into  $A=\{a_{i}|0 \leq i \leq 7\}$ and $B=\{b_{i}|0 \leq i \leq 7\}$. First we decompose $K_{8,8}$ into cycles of length $8$. Then we pick two or three cycles together and decompose them into paths of length $8$. Let $ C^{1}=(a_{1}  b_{1}  a_{7}  b_{7}  a_{5}  b_{5}  a_{3}  b_{3}),\ C^{2}=(a_{1}  b_{0}  a_{7}  b_{6}  a_{5}  b_{4}  a_{3}  b_{2}),\ C^{3}=(a_{1}  b_{4}  a_{7}  b_{2}  a_{5}  b_{0} a_{3}  b_{6})$ and $C^{4}=(a_{1}  b_{5}  a_{7}  b_{3}  a_{5}  b_{1}  a_{3}  b_{7})$. Consider the permutation $\rho=(a_{0}a_{1}a_{2}a_{3}a_{4} a_{5}a_{6}a_{7})(b_{0}b_{1}b_{2}b_{3}b_{4}b_{5}b_{6}b_{7})$. Then $\{\rho( C^{i})| 1 \leq i \leq 4\}$ provides desired $C_{8}$-decomposition. Now, $(8; p, q)$-decomposition of $K_{8,8}$ is as follows:\\
	$D_{1}$: Decomposition of $C^{1}\cup C^{2}$ into $2$ copies of $ P_{8}$ given by:  $P^{1}= b_{0}  a_{7}  b_{7}  a_{5}  b_{5}  a_{3}\\ b_{3}  a_{1}  b_{1}$ and $P^{2}=b_{0}  a_{1}  b_{2}  a_{3}  b_{4}  a_{5}  b_{6}  a_{7}  b_{1}$.\\
	$D_{2}$: Similarly $\rho(P^{1})$ and $\rho(P^{2})$ are the $8$-paths obtained from $\rho( C^{1})\cup \rho( C^{2})$.\\
	$D_{3}$: Decomposition of $C^{3}\cup C^{4}$ into $2$ copies of $ P_{8}$ given by: $P^{3}= b_{7}  a_{1}  b_{6}  a_{3}  b_{0} \\ a_{5}  b_{2}  a_{7}  b_{4}$ and $P^{4}=b_{4}  a_{1}  b_{5}  a_{7}  b_{3}  a_{5}  b_{1}  a_{3}  b_{7}$.\\
	$D_{4}$: Similarly $\rho(P^{3})$ and $\rho(P^{4})$ are the $8$-paths obtained from $\rho( C^{3})\cup \rho( C^{4})$.\\
	$D_{5}$: Consider $C^{1}\cup C^{2}\cup C^{3}$ and decompose them into $3$ copies of $ P_{8}$ as follows:  $P^{2}=b_{0}  a_{1}  b_{2}  a_{3}  b_{4}  a_{5}  b_{6}  a_{7}  b_{1}, \ P^{5}=b_{1}  a_{1}  b_{4}  a_{7}  b_{2}  a_{5}  b_{0}  a_{3}  b_{6}$ and $P^{6}=b_{6}  a_{1}  b_{3}  a_{3}  b_{5}  a_{5}  b_{7} a_{7}  b_{0}$. \\
	$D_{6}$: Similarly $\rho(P^{2}),\ \rho(P^{5})$ and $\rho(P^{6})$ are the $8$-paths obtained from $\rho( C^{1})\cup \rho( C^{2})\cup \rho( C^{3})$.\\
	Here $(p,q) \in \{(0,8),(2,6),(3,5),(4,4),(5,3),(6,2),(7,1),(8,0)\}$. For $1 \leq i \leq 6$, combinations of $D_{i}$ provide the required $(8; p, q)$-decomposition.\hfill$\Box$
	
	\begin{lemma}
		The graph $P_{2} \times K_{5}$ has a $(8; p, q)$-decomposition.
		\label{LEM4C8P8}
	\end{lemma}
	\noindent{\bf Proof.}
	Let the vertices of the graph $P_{2} \times K_{5}$ be partitioned into $A=\{a_{i}|0 \leq i \leq 4\},\ B=\{b_{i}|0 \leq i \leq 4\}$ and $C=\{c_{i}|0 \leq i \leq 4\}$. The $C_{8}$-decomposition of $P_{2} \times K_{5}$ given by $\{\rho^{i}(C) | 0 \leq i \leq 4 \}$, where $C=(a_{0}b_{4}a_{2}b_{0}c_{1}b_{3}c_{0}b_{1})$ and $\rho=(a_{0}a_{1}a_{2}a_{3}a_{4})(b_{0}b_{1}b_{2}b_{3}b_{4})(c_{0}c_{1}c_{2}c_{3}c_{4})$. The required $(8; p, q)$-decomposition of $P_{2} \times K_{5}$ is as follows:\\
	$D_{1}$: By decomposing $C\cup \rho(C)$ we obtain $2$ copies of $ P_{8}$ as follows: $P^{1}=a_{0} b_{4} a_{2} b_{0} c_{1} b_{3} c_{0} b_{1} c_{2}$ and $P^{2}=a_{0} b_{1} a_{3} b_{0} a_{1} b_{2} c_{1} b_{4} c_{2}$.\\
	$D_{2}$: $\rho^{3}(P^{1})$ and $\rho^{3}(P^{2})$ provide $P_{8}$ decomposition of $\rho^{3}(C) \cup \rho^{4}(C)$. \\
	$D_{3}$: Decomposition of $ C \cup \rho(C)\cup \rho^{2}(C)$ into $3$ copies of $ P_{8}$ given by: $P^{1}=a_{0} b_{4} a_{2} b_{0} c_{1} b_{3} c_{0} b_{1} c_{2},\ P^{3}=a_{0}  b_{1}  a_{2}  b_{3}  c_{2}  b_{0}  c_{3} b_{2} a_{4}$ and $ P^{4}= a_{4}  b_{1}  a_{3}  b_{0}  a_{1}  b_{2}  c_{1}  b_{4} c_{2}$.\\
	Here $(p,q) \in \{(0,5),(2,3),(3,2),(4,1),(5,0)\}$. For $1 \leq i \leq 3$, combinations of $D_{i}$ provide the required $(8; p, q)$-decomposition.\hfill$\Box$
	
	\begin{lemma}
		There exists a $(8; p, q)$-decomposition of $K_{3} \times P_{4}$.
		\label{LEM5C8P8}
	\end{lemma}
	\noindent{\bf Proof.}
	Partition the vertices of the graph $K_{3} \times P_{4}$ into $A=\{a_{i}|0 \leq i \leq 4\},\ B=\{b_{i}|0 \leq i \leq 4\}$ and $C=\{c_{i}|0 \leq i \leq 4\}$. The decomposition of $K_{3} \times P_{4}$ into $8$-cycles given by: $C^{1}=(c_{0}a_{1}b_{2}a_{3}c_{4}b_{3}c_{2}b_{1}),\ C^{2}=(c_{1}b_{0}a_{1}c_{2}a_{3}b_{4}c_{3}a_{2})$ and $C^{3}=(c_{1}a_{0}b_{1}a_{2}b_{3}a_{4}c_{3}b_{2})$. The required $(8; p, q)$-decomposition of $K_{3} \times P_{4}$ is as follows:\\
	For $(p,q)=(2,1)$: Decompose $C^{1}\cup C^{3}$ into $2$ copies of $ P_{8}$, say $ c_{1}  b_{2}  a_{3}  c_{4}  b_{3}  c_{2}  b_{1} \\ c_{0} a_{1}$ and $ a_{1}  b_{2}  c_{3}  a_{4}  b_{3}  a_{2} b_{1}  a_{0} c_{1}$.\\
	For $(p,q)=(3,0)$: By decomposing $C^{1}\cup C^{2}\cup C^{3}$ we obtain $3$ copies of $ P_{8}$ as follows: $ c_{1}  b_{0}  a_{1}  c_{0}  b_{1}  c_{2}  a_{3}  b_{4} c_{3},\ c_{3}  a_{2}  c_{1}  b_{2}  a_{3}  c_{4}  b_{3}  c_{2} a_{1}$ and $a_{1}  b_{2}  c_{3}  a_{4}  b_{3}  a_{2}  b_{1}  a_{0} c_{1}$.\hfill$\Box$ 
	
	\begin{lemma}
		The graph $ P_{2} \times K_{8}$ has a $(8; p, q)$-decomposition.
		\label{LEM6C8P8}
	\end{lemma}
	\noindent{\bf Proof.}
	Partition the vertices of the graph $P_{2} \times K_{8}$ be $A=\{a_{i}|0 \leq i \leq 7\},\ B=\{b_{i}|0 \leq i \leq 7\}$ and $C=\{c_{i}|0 \leq i \leq 7\}$. Consider $P=a_{3}b_{1}c_{0}b_{5}a_{1}b_{3}c_{4}b_{7}c_{1}$ and the permutation $\rho=(a_{0}a_{1}a_{2}a_{3}a_{4}a_{5}a_{6}a_{7})(b_{0}b_{1}b_{2}b_{3}b_{4}b_{5}\\b_{6}b_{7})(c_{0}c_{1}c_{2}c_{3}c_{4}c_{5}c_{6}c_{7})$. Then $\{\rho^{i}(P) | 0 \leq i \leq 7 \}$ provides $8$ copies of $ P_{8}$. The $C_{8}$-decomposition of remaining edges of $ P_{2} \times K_{8}$ given by: $\{\rho^{i}(C^{j}) | i=0,1$ and $ 1 \leq j \leq 3 \}$, where $C^{1}=(a_{0}b_{1}a_{2}b_{3}a_{4}b_{5}a_{6}b_{7}),\ C^{2}=(a_{0}b_{3}a_{6}b_{1}a_{4}b_{7}a_{2}b_{5}) $ and $C^{3}=(b_{0}c_{4}b_{6}c_{2}b_{4}c_{0}b_{2}c_{6})$. The required $(8; p, q)$-decomposition of $ P_{2} \times K_{8}$ is as follows:\\
	$D_{1}$: Decomposition of $C^{1}\cup P$ into $2$ copies of $ P_{8}$ given by: $P^{1}= a_{3}  b_{1}  a_{0} b_{7}  a_{6}  b_{5} \\ a_{4} b_{3} a_{2}$ and $P^{2}= a_{2}  b_{1}  c_{0}  b_{5}  a_{1}  b_{3}  c_{4}  b_{7} c_{1}$. \\
	$D_{2}$: By decomposing $C^{1}\cup \rho(C^{3})$, we obtain $2 8$-paths, say $P^{3}= a_{0}  b_{1}  c_{7} b_{3}  c_{1}  b_{5}  c_{3} \\ b_{7} c_{5}$ and $P^{4}= c_{5}  b_{1}  a_{2}  b_{3}  a_{4}  b_{5}  a_{6}  b_{7} a_{0}$. \\
	$D_{3}$: Consider $\rho (C^{1})\cup C^{3}$, then $\rho(P^{3})$ and $\rho(P^{4})$ provides a $P_{8}$-decomposition.\\
	$D_{4}$: Decomposition of $C^{1}\cup C^{2}\cup \rho(C^{3})$ into $3$ copies of $ P_{8}$ given by: $ P^{5}= c_{7}  b_{1}  a_{6}  b_{3}  a_{0}  b_{5}  a_{2}  b_{7} a_{4}, \ P^{6}= c_{7}  b_{3}  a_{4}  b_{5}  c_{3}  b_{7}  c_{5}  b_{1} a_{2}$ and $ P^{7}= a_{4}  b_{1}  a_{0} b_{7}  a_{6}  b_{5}  c_{1}  b_{3} a_{2}$. \\
	$D_{5}$: Similarly $\rho(P^{5}),\ \rho(P^{6})$ and $\rho(P^{7})$ establish a $P_{8}$-decomposition of  $\rho (C^{1})\cup \rho (C^{2})\cup C^{3}$. \\
	Here $(p,q) \in \{(8,6),(9,5),(10,4),(11,3),(12,2),(13,1),(14,0)\}$. For $1 \leq i \leq 5$, combinations of $D_{i}$ provide the required $(8; p, q)$-decomposition.\hfill$\Box$
	
	\begin{lemma}
		The graph $C_{4} \times K_{4}$ has a $(8; p, q)$-decomposition.
		\label{LEMC_{4}K_{4}C8P8}
	\end{lemma}
	\noindent{\bf Proof.}
	Let the vertices of the graph $C_{4} \times K_{4}$ be partitioned into $A=\{a_{i}|0 \leq i \leq 3\},\ B=\{b_{i}|0 \leq i \leq 3\},\ C=\{c_{i}|0 \leq i \leq 3\}$ and $D=\{d_{i}|0 \leq i \leq 4\}$. The decomposition of $C_{4} \times K_{4}$ into $8$-cycles given by: $C^{1}=(a_{0}b_{1}c_{0}d_{1}c_{3}d_{2}c_{1}b_{3}),\ C^{2}=(a_{0}d_{1}c_{2}b_{0}a_{3}b_{2}c_{0}d_{2}),\ C^{3}=(d_{0}a_{1}d_{2}a_{3}d_{1}a_{2}d_{3}c_{2}),\ C^{4}=(b_{0}a_{2}b_{3}c_{2}b_{1}a_{3}d_{0}c_{3}),\ C^{5}=(d_{0}c_{1}d_{3}a_{1}b_{2}c_{3}b_{1}a_{2})$ and $C^{6}=(a_{0}b_{2}c_{1}b_{0}a_{1}b_{3}c_{0}d_{3})$. The required $(8; p, q)$-decomposition of $C_{4} \times K_{4}$ is as follows:\\
	$D_{1}$: Consider $C^{1}\cup C^{2}$ and decompose them into $2$ copies of $ P_{8}$ given by: $P^{1}=c_{2}  d_{1} c_{0}  b_{1}  a_{0}  b_{3}  c_{1} d_{2} c_{3}$ and $P^{2}=c_{3}  d_{1} a_{0}  d_{2}  c_{0}  b_{2}  a_{3} b_{0} c_{2}$.\\
	$D_{2}$: By decomposing $C^{4}\cup C^{5}$ we obtain $P_{8}$ as follows: $P^{3}=b_{2}  c_{3} d_{0}  a_{3}  b_{1}  c_{2}  b_{3} a_{2} \\b_{0}$ and $P^{4}=b_{0}  c_{3} b_{1}  a_{2}  d_{0}  c_{1}  d_{3} a_{1} b_{2}$.\\
	$D_{3}$: Decomposition of $C^{3}\cup C^{6}$ into $2$ copies of $ P_{8}$ given by: $P^{5}=b_{0}  a_{1} d_{0}  c_{2}  d_{3}  a_{2}  \\d_{1}  a_{3} d_{2}$ and $P^{6}=d_{2}  a_{1} b_{3}  c_{0}  d_{3}  a_{0}  b_{2} c_{1} b_{0}$.\\
	$D_{4}$: Decompose $C^{1}\cup C^{2}\cup C^{6}$ into $P_{8}$, say $P^{2}=c_{3}  d_{1} a_{0}  d_{2}  c_{0}  b_{2}  a_{3} b_{0} c_{2}, P^{7}=c_{3}  d_{2} c_{1}  b_{2}  a_{0}  d_{3}  c_{0} b_{3} a_{1}$ and $P^{8}=a_{1}  b_{0} c_{1}  b_{3}  a_{0}  b_{1}  c_{0} d_{1} c_{2}$.\\
	Here $(p,q) \in \{(0,6),(2,4),(3,3),(4,2),(5,1),(6,0)\}$. For $1 \leq i \leq 4$, combinations of $D_{i}$ provide the required $(8; p, q)$-decomposition.\hfill$\Box$
	
	\begin{lemma}
		The graph $K_{2} \times K_{9}$ has a $(8; p, q)$-decomposition.
		\label{LEM7C8P8}
	\end{lemma}
	\noindent{\bf Proof.}
	Partition the vertices of the graph $K_{2} \times K_{9}$ into $A=\{a_{i}|0 \leq i \leq 8\}$ and $B=\{b_{i}|0 \leq i \leq 8\}$. Consider $C^{1}=(a_{0}b_{2}a_{1}b_{5}a_{2}b_{0}a_{3}b_{8})$ and the permutation $\rho=(a_{0}a_{1}a_{2}a_{3}a_{4}a_{5}a_{6}a_{7}a_{8})(b_{0}b_{1}b_{2}b_{3}b_{4}b_{5}b_{6}b_{7}b_{8})(c_{0}c_{1}c_{2}c_{3}c_{4}c_{5}c_{6}c_{7}c_{8})$. \\Then $\{\rho^{i} (C^{1}) | 1\leq i \leq 8 \}$ provides $ 8$-cycles, namely $C^{2},C^{3},C^{4},C^{5},C^{6},C^{7},C^{8}$ and $C^{9}$. The required $(8; p, q)$-decomposition of $K_{2} \times K_{9}$ is as follows:\\
	$D_{1}$: Decomposition of $C^{1}\cup C^{2}$ into $2$ copies of $ P_{8}$ given by: $P^{1}=b_{3}  a_{2} b_{6}  a_{3}  b_{1}  a_{4} \\ b_{0} a_{1} b_{2}$ and $P^{2}=b_{2}  a_{0} b_{8}  a_{3}  b_{0}  a_{2}  b_{5} a_{1} b_{3}$.\\
	$D_{2}$: Consider $C^{3}\cup C^{4}$, then $\rho^{2}(P^{1})$ and $\rho^{2}(P^{2})$ provide $P_{8}$-decomposition.\\
	$D_{3}$: Consider $C^{5}\cup C^{6}$, then $\rho^{4}(P^{1})$ and $\rho^{4}(P^{2})$ provide $P_{8}$-decomposition.\\
	$D_{4}$: Consider $C^{7}\cup C^{8}$, then $\rho^{6}(P^{1})$ and $\rho^{6}(P^{2})$ provide $P_{8}$-decomposition.\\
	$D_{5}$: Decompose $C^{1}\cup C^{2}\cup C^{9}$ into $P_{8}$ as follows: $P^{2}=b_{2}  a_{0} b_{8}  a_{3}  b_{0}  a_{2}  b_{5} a_{1} b_{3}, \\P^{3}=b_{3}  a_{2} b_{1}  a_{5}  b_{2}  a_{4}  b_{7} a_{3} b_{4}$ and $P^{4}=b_{4}  a_{2} b_{6}  a_{3}  b_{1}  a_{4}  b_{0} a_{1} b_{2}$.\\
	Here $(p,q) \in \{(0,9),(2,7),(3,6),(4,5),(5,4),(6,3),(7,2),(8,1),(9,0)\}$. For $1 \leq i \leq 5$, combinations of $D_{i}$ provide the required $(8; p, q)$-decomposition.\hfill$\Box$
	
	\begin{lemma}
		The graph $C_{4} \times K_{6}$ admits a $(8; p, q)$-decomposition.
		\label{LEMC_{4}K_{6}C8P8}
	\end{lemma}
	\noindent{\bf Proof.}
	Let the vertices of the graph $C_{4} \times K_{6}$ be partitiones into $A=\{a_{i}|0 \leq i \leq 5\},\ B=\{b_{i}|0 \leq i \leq 5\},\ C=\{c_{i}|0 \leq i \leq 5\}$ and $D=\{d_{i}|0 \leq i \leq 5\}$. The decomposition of $C_{4} \times K_{6}$ into $8$-cycles given by: $C^{1}=(b_{0}c_{1}d_{2}c_{3}b_{5}c_{0}d_{5}c_{4}),\ C^{2}=(a_{2}b_{0}c_{3}d_{5}a_{0}d_{1}c_{4}b_{1}),\ C^{3}=(a_{1}d_{3}c_{4}b_{0}a_{4}b_{5}a_{2}d_{5 })$ and the permutation $\rho=(a_{0}a_{1}a_{2}a_{3}a_{4}a_{5})(b_{0}b_{1}b_{2}b_{3}b_{4}b_{5})(c_{0}c_{1}c_{2}c_{3}c_{4}c_{5})(d_{0}d_{1}d_{2}d_{3}d_{4} \\ d_{5})$. Then $\{\rho^{i} (C^{j}) | 0\leq i \leq 4,\ j=1,2,3 \}$ provides a $C_{8}$ decomposition of $C_{4} \times K_{6}$. The required $(8; p, q)$-decomposition of $C_{4} \times K_{6}$ is as follows: \\
	$D_{1}$: Decomposition of $C^{1}\cup C^{2}$ into $2$ copies of $ P_{8}$ is given by $P^{1}=a_{2}  b_{0} c_{4}  d_{5}  c_{0} \\ b_{5}  c_{3} d_{2} c_{1}$ and $P^{2}=c_{1}  b_{0} c_{3}  d_{5}  a_{0}  d_{1}  c_{4} b_{1} a_{2}$.\\
	$D_{2}$: For $1\leq i \leq 4$, $\rho^{i}(P^{1})$ and $\rho^{i}(P^{2})$ provides a $P_{8}$-decomposition of $\rho^{i}(C^{1})\cup \rho^{i}(C^{2})$.\\
	$D_{3}$: Decompose $C^{1}\cup C^{2}\cup C^{3}$ into $P_{8}$ as follows: $P^{2}=c_{1}  b_{0} c_{3}  d_{5}  a_{0}  d_{1}  c_{4} b_{1} a_{2},\\ P^{3}=c_{1}  d_{2} c_{3}  b_{5}  a_{2}  d_{5}  c_{4} b_{0} a_{4}$ and $P^{4}=a_{4}  b_{5} c_{0}  d_{5}  a_{1}  d_{3}  c_{4} b_{0} a_{2}$.\\
	$D_{4}$: For $1\leq i \leq 4$, $\rho^{i}(P^{1}),\ \rho^{i}(P^{2})$ and $\rho^{i}(P^{3})$ provides a $P_{8}$-decomposition of $\rho^{i}(C^{1})\cup \rho^{i}(C^{2}) \cup \rho^{i}(C^{3})$.\\
	Here $(p,q) \in \{(0,15),(2,13),(3,12),(4,11),(5,10),(6,9),(7,8),(8,7),(9,6),\\(10,5),(11,4),(12,3),(13,2),(14,1),(15,0)\}$. For $1 \leq i \leq 4$, combinations of $D_{i}$ provide the required $(8; p, q)$-decomposition.\hfill$\Box$
	
	\begin{lemma}
		The graph $C_{4} \times K_{3}$ admits a $(8; p, q)$-decomposition.
		\label{LEMC_{4}K_{3}C8P8}
	\end{lemma}
	\noindent{\bf Proof.}
	Let the vertices of the graph $C_{4} \times K_{3}$ be partitioned into $A=\{a_{i}|0 \leq i \leq 2\},\ B=\{b_{i}|0 \leq i \leq 2\},\ C=\{c_{i}|0 \leq i \leq 2\}$ and $D=\{d_{i}|0 \leq i \leq 2\}$. The decomposition of $C_{4} \times K_{3}$ into $8$-cycles given by: $C^{1}=(a_{1}d_{0}c_{2}b_{1}a_{2}b_{0}c_{1}d_{2}),\ C^{2}=(a_{0}b_{1}c_{0}d_{1}a_{2}d_{0}c_{1}b_{2})$ and $C^{3}=(a_{1}b_{0}c_{2}d_{1}a_{0}d_{2}c_{0}b_{2})$. The required $(8; p, q)$-decomposition of $C_{4} \times K_{3}$ is as follows:\\
	For $(p,q)=(2,1)$: Decomposition of $C^{1}\cup C^{2}$ into $2$ copies of $ P_{8}$ given by: $P^{1}=a_{0}  b_{1} a_{2}  b_{0}  c_{1}  d_{2}  a_{1} d_{0} c_{2}$ and $P^{2}=c_{2}  b_{1} c_{0}  d_{1}  a_{2}  d_{0}  c_{1} b_{2} a_{0}$.\\
	For $(p,q)=(3,0)$: Decompose $C^{1}\cup C^{2}\cup C^{3}$ into $P_{8}$ as follows: $P^{2}=c_{2}  b_{1} c_{0}  d_{1}  a_{2}  d_{0}  c_{1} b_{2} a_{0}, P^{3}=a_{2}  b_{1} a_{0}  d_{2}  c_{0}  b_{2}  a_{1} b_{0} c_{2}$ and $P^{4}=a_{2}  b_{0} c_{1}  d_{2}  a_{1}  d_{0}  c_{2} d_{1} a_{0}$.\hfill$\Box$
	
	\begin{lemma}
		The graph $(K_{3} \times K_{3})(4)$ admits a $(8; p, q)$-decomposition.
		\label{LEM4K_{3}K_{3}C8P8}
	\end{lemma}
	\noindent{\bf Proof.}
	Let the vertices of the graph $K_{3} \times K_{3}$ be partitioned into $A=\{a_{i}|0 \leq i \leq 2\},\ B=\{b_{i}|0 \leq i \leq 2\}$ and $C=\{c_{i}|0 \leq i \leq 2\}$. The decomposition of $(K_{3} \times K_{3})(4)$ into $8$-cycles given by: $C^{1}=(a_{0}c_{1}b_{0}a_{1}b_{2}c_{0}b_{1}c_{2}),\ C^{2}=(a_{0}c_{1}b_{2}a_{1}c_{2}b_{0}a_{2}b_{1}),\ C^{3}=(a_{0}b_{1}c_{2}a_{1}c_{0}a_{2}c_{1}b_{2}),\ C^{4}=(a_{0}c_{2}b_{1}c_{0}a_{2}b_{0}c_{1}b_{2}),\ C^{5}=(a_{0}c_{1}a_{2}b_{1}c_{0}a_{1}b_{0}c_{2}),\ C^{6}=(a_{0}b_{1}a_{2}b_{0}c_{2}a_{1}c_{0}b_{2}),\ C^{7}=(a_{0}c_{1}a_{2}c_{0}b_{2}a_{1}b_{0}c_{2}),\ C^{8}=(a_{0}b_{1}c_{0}a_{2}c_{1}b_{0}a_{1}b_{2})$ and $C^{9}=(a_{1}c_{2}b_{1}a_{2}b_{0}c_{1}b_{2}c_{0}) $. The required $(8; p, q)$-decomposition of $(K_{3} \times K_{3})(4)$ is as follows:\\
	$D_{1}$: Decomposition of $C^{1}\cup C^{2}$ into $2$ copies of $ P_{8}$ given by: $P^{1}=a_{2}  b_{1} c_{2}  a_{0}  c_{1}  b_{0} \\ a_{1} b_{2} c_{0}$ and $P^{2}=c_{0}  b_{1} a_{0}  c_{1}  b_{2}  a_{1}  c_{2} b_{0} a_{2}$.\\
	$D_{2}$: Decomposition of $C^{3}\cup C^{8}$ into $2$ copies of $ P_{8}$ given by: $P^{3}=c_{2}  b_{1} a_{0}  b_{2}  c_{1}  a_{2} \\ c_{0} a_{1} b_{0}$ and $P^{4}=c_{2}  a_{1} b_{2}  a_{0}  b_{1}  c_{0}  a_{2} c_{1} b_{0}$.\\
	$D_{3}$: Decomposition of $C^{4}\cup C^{5}$ into $2$ copies of $ P_{8}$ given by: $P^{5}=a_{1}  c_{0} b_{1}  a_{2}  b_{0}  c_{2} \\ a_{0} c_{1} b_{2}$ and $P^{6}=b_{2}  a_{0} c_{2}  b_{1}  c_{0}  a_{2}  c_{1} b_{0} a_{1}$.\\
	$D_{4}$: Decomposition of $C^{6}\cup C^{7}$ into $2$ copies of $ P_{8}$ given by: $P^{7}=c_{1}  a_{2} b_{0}  c_{2}  a_{1}  c_{0} \\ b_{2} a_{0} b_{1}$ and $P^{8}=b_{1}  a_{2} c_{0}  b_{2}  a_{1}  b_{0}  c_{2} a_{0} c_{1}$.\\
	$D_{5}$: Decompose $C^{1}\cup C^{2}\cup C^{9}$ into $P_{8}$ as follows: $P^{1}=a_{2}  b_{1} c_{2}  a_{0}  c_{1}  b_{0}  a_{1} b_{2} c_{0},\\ P^{9}=a_{0}  b_{1} c_{2}  a_{1}  c_{0}  b_{2}  c_{1} b_{0} a_{2}$ and $P^{10}=a_{0}  c_{1} b_{2}  a_{1}  c_{2}  b_{0}  a_{2} b_{1} c_{0}$.\\
	Here $(p,q) \in \{(0,9),(2,7),(3,6),(4,5),(5,4),(6,3),(7,2),(8,1),(9,0)\}$. For $1 \leq i \leq 5$, combinations of $D_{i}$ provide the required $(8; p, q)$-decomposition.\hfill$\Box$

	\section{$(8; p, q)$-decomposition of $ (K_{n} \times K_{m})(\lambda)$}
	In this section, we prove the necessary and sufficient conditions for $(8; p, q)$-decomposition of $(K_{n} \times K_{m})(\lambda)$.\\ 
	\begin{theorem}
		For $mn \geq 9$, $K_{n} \times K_{m}(\lambda)$ admits a $(8; p, q)$-decomposition if and only if $8|\frac{\lambda mn(m-1)(n-1)}{2}$.
	\end{theorem}
	\noindent{\bf Proof.}
	The necessary conditions are obvious and we prove the sufficiency as follows:\\
	
	\noindent{\bf Case 1:} $\lambda=1.$\\
	\noindent{\bf Subcase 1.1:} $n \equiv 0($mod $4)$ and $m \equiv 0($mod $4)$.\\
	\noindent Here $K_{m}$ and $K_{n}$ is decomposed into subgraphs $K_{4}$ and $K_{4,4}$. Hence $K_{n} \times K_{m}$ is isomorphic to subgraphs $(K_{4} \times K_{4}),\ (K_{4,4} \times K_{4,4})$ and $(K_{4,4} \times K_{4})$. Clearly, $K_{4,4} \times K_{4,4}$ and $K_{4,4} \times K_{4}$ can be expressed in terms of $ P_{2} \circ \overline{K}_{4}$ which admits a $(8; p, q)$-decomposition by Lemma \ref{LEM2C8P8}. By Lemma \ref{LEM1C8P8}, $K_{4} \times K_{4}$ has a $(8; p, q)$-decomposition. Hence $K_{n} \times K_{m}$ admits a $(8; p, q)$-decomposition. \\ 
	
	\noindent{\bf Subcase 1.2:} $n \equiv 0($mod $4)$ and $m \equiv 1($mod $4)$.\\
	\noindent Here $K_{m}\cong \left(\frac{m-1}{4}\right) K_{5} \oplus \left(\frac{(m-1)(m-5)}{32}\right)K_{4,4}$ and $K_{n}\cong \left(\frac{n}{4}\right)K_{4} \oplus \left(\frac{n^{2}-4n}{32}\right)K_{4,4}$. By Theorem \ref{THM1C8P8}, $P_{2}|K_{4}$, so we rewrite $K_{n} \times K_{m} \cong \left(\frac{3n(m-1)}{16}\right)(P_{2} \times K_{5}) \oplus \left(\frac{3n(m-1)(m-5)}{128}\right)(P_{2} \times K_{4,4}) \oplus \left(\frac{n(n-4)(m-1)}{128}\right)(K_{5} \times K_{4,4}) \oplus \left(\frac{n(n-4)(m-1)(m-5)}{1024}\right) \\ (K_{4,4} \times K_{4,4})$. As $P_{2} \times K_{4,4},\ K_{5} \times K_{4,4}$ and $K_{4,4} \times K_{4,4}$ can be expressed in terms of $P_{2} \circ \overline{K}_{4}$ which admits a $(8; p, q)$-decomposition by Lemma \ref{LEM2C8P8}. By Lemma \ref{LEM4C8P8}, $P_{2} \times K_{5}$ has a $(8; p, q)$-decomposition. \\
	
	\noindent{\bf Subcase 1.3:} $n \equiv 1($mod $4)$ and $m \equiv 1($mod $4)$.\\
	\noindent Consider $K_{m}\cong \left(\frac{m-1}{4}\right)  K_{5} \oplus \left(\frac{(m-1)(m-5)}{32}\right) K_{4,4} $ and $K_{n}\cong \left(\frac{n-1}{4}\right)  K_{5} \oplus \left(\frac{(n-1)(n-5)}{32}\right) K_{4,4} $. By Theorem \ref{THM1C8P8}, $P_{2}|K_{5}$, then we have $K_{n} \times K_{m} \cong \\ \left(\frac{5(n-1)(m-1)}{16}\right) (P_{2} \times K_{5}) \oplus \left(\frac{5(n-1)(m-1)(m-5)}{128}\right)  (P_{2} \times K_{4,4}) \oplus \left(\frac{(n-1)(n-5)(m-1)}{128}\right) \\ (K_{5} \times K_{4,4}) \oplus \left(\frac{(n-1)(n-5)(m-1)(m-5)}{1024}\right)  (K_{4,4} \times K_{4,4})$. Clearly $P_{2} \times K_{4,4},\ K_{5} \times K_{4,4}$ and $K_{4,4} \times K_{4,4}$ can be expressed in terms of $P_{2} \circ \overline{K}_{4}$ which admits a $(8; p, q)$-decomposition by Lemma \ref{LEM2C8P8}. By Lemma \ref{LEM4C8P8}, $P_{2} \times K_{5}$ has a $(8; p, q)$-decomposition.\\
	
	\noindent{\bf Subcase 1.4:} $n \equiv 0($mod $8)$ and $m \equiv 0($mod $2)$.\\
	\noindent{\bf Subcase 1.4.1:} $n \equiv 0($mod $8)$ and $m \equiv 0($mod $4)$.\\
	\noindent The proof of this case follows from Subcase 1.1.\\ 
	
	\noindent {\bf Subcase 1.4.2:} $n \equiv 0($mod $8)$ and $m \equiv 2($mod $4)$.\\
	\noindent Here $K_{n}\cong \left(\frac{n}{8}\right) K_{8} \oplus \left(\frac{n^{2}-8n}{32}\right) K_{4,4}$ and $K_{m}\cong K_{6} \oplus \left(\frac{m-6}{4}\right) (K_{6}-e) \oplus \left(\frac{(m-2)(m-6)}{32}\right) K_{4,4}$. Now, $K_{6} \cong (K_{4}-e) \oplus K_{4} \oplus C_{4}$ and $(K_{6}-e) \cong 2 (K_{4}-e) \oplus C_{4}$. Hence $K_{n} \times K_{m} \cong \left(\frac{n}{8}\right)  (K_{8} \times K_{4}) \oplus \left(\frac{n(m-6)(m-2)}{256}\right)  (K_{8} \times K_{4,4}) \oplus \left(\frac{n(m-4)}{16}\right)  (K_{8} \times (K_{4}-e)) \oplus \left(\frac{n(m-2)}{32}\right)  (K_{8} \times C_{4}) \oplus \left(\frac{n(n-8)}{32}\right)  (K_{4,4} \times K_{4})\oplus \left(\frac{n(n-8)(m-6)(m-2)}{1024}\right)  (K_{4,4} \times K_{4,4})\oplus \left(\frac{n(n-8)(m-4)}{64}\right)  (K_{4,4} \times (K_{4}-e))\oplus \left(\frac{n(n-8)(m-2)}{128}\right)  (K_{4,4} \times C_{4})$. Clearly, $K_{8} \times K_{4,4}, K_{4,4} \times K_{4} ,K_{4,4} \times K_{4,4}$ and $K_{4,4} \times C_{4}$ can be expressed in terms of $P_{2} \circ \overline{K}_{4}$ which admits a $(8; p, q)$-decomposition by Lemma \ref{LEM2C8P8}. Let $K_{8} \times K_{4} \cong (2 K_{4} \oplus K_{4,4}) \times K_{4} $. By Lemma \ref{LEM1C8P8}, $K_{4}\times K_{4}$ admits a $(8; p, q)$-decomposition and $K_{4,4}\times K_{4}$ can be expressed in terms of $P_{2} \circ \overline{K}_{4}$ which admits a $(8; p, q)$-decomposition by Lemma \ref{LEM2C8P8}. Let $K_{8} \times C_{4} \cong (2 K_{4} \oplus K_{4,4})  \times C_{4}$. By Lemma \ref{LEMC_{4}K_{4}C8P8}, $K_{4}\times C_{4}$ admits a $(8; p, q)$-decomposition and $K_{4,4}\times C_{4}$ can be expressed in terms of $P_{2} \circ \overline{K}_{4}$ which admits a $(8; p, q)$-decomposition by Lemma \ref{LEM2C8P8}. Since $(K_{4}-e)$ can be decomposed into $P_{2}$ and $K_{3}$ and $P_{4}|K_{8}$  by Theorem \ref{THM1C8P8}, we write $K_{8} \times (K_{4}-e) \cong (K_{8}\times P_{2}) \oplus (K_{8}\times K_{3}) \cong(K_{8}\times P_{2}) \oplus (7P_{4}\times K_{3})$. Hence $K_{8} \times (K_{4}-e)$ has a $(8; p, q)$-decomposition by Lemma \ref{LEM5C8P8} and \ref{LEM6C8P8}. Let $K_{4,4} \times (K_{4}-e) \cong (K_{4,4}\times P_{2}) \oplus (K_{4,4}\times K_{3}) \cong 2 (P_{2} \circ \overline{K}_{4}) \oplus 4 (P_{4}\times K_{3}) $. Hence $K_{4,4} \times (K_{4}-e)$ has a $(8; p, q)$-decomposition by Lemmas \ref{LEM2C8P8} and \ref{LEM5C8P8}.\\
	
	\noindent{\bf Subcase  1.5:} $n \equiv 0($mod $8)$ and $m \equiv 1($mod $2)$.\\
	\noindent{\bf Subcase 1.5.1:} $n \equiv 0($mod $8)$ and $  m \equiv 1$ and $ 3($mod $ 6) $.\\
	\noindent Let $K_{n} \cong \left(\frac{n}{8}\right) K_{8} \oplus \left(\frac{n^{2}-8n}{128}\right) K_{8,8} \cong \left(\frac{7n}{8}\right) P_{4} \oplus \left(\frac{n^{2}-8n}{128}\right) K_{8,8}$. By Theorem \ref{THM2C8P8}, $K_{3}|K_{m}$. By Lemma \ref{LEM5C8P8}, $P_{4}\times K_{3}$ has a $(8; p, q)$-decomposition. By Lemmas \ref{LEM0C8P8} and \ref{LEM3C8P8}, $ K_{8,8}\times K_{3} $ has a $(8; p, q)$-decomposition. \\
	
	\noindent{\bf Subcase 1.5.2:} $n \equiv 0($mod $8)$ and $  m \equiv 5($mod $ 6 ) $.\\
	\noindent Here $K_{n}\cong \left(\frac{n}{8}\right) K_{8} \oplus \left(\frac{n^{2}-8n}{128}\right) K_{8,8}$. By Theorem \ref{THM2C8P8}, $K_{m} \cong \left(\frac{m(m -1) - 20}{6}\right) K_{3} \oplus K_{5} $. Then $K_{n} \times K_{m} \cong \left(\frac{n(m(m -1) - 20)}{48}\right) (K_{8}\times K_{3}) \oplus \left(\frac{n}{8}\right) (K_{8}\times K_{5})\oplus \left(\frac{(n^{2}-8n)(m(m -1) - 20)}{768}\right) (K_{8,8}\times K_{3}) \oplus \left(\frac{n(n-8)}{128}\right) (K_{8,8}\times K_{5}))$. Now $K_{8}\times K_{3}$ can be expressed in terms of $P_{4}\times K_{3}$, which has a $(8; p, q)$-decomposition by Lemma \ref{LEM5C8P8}. By Theorem \ref{THM1C8P8}, $K_{8}\times K_{5}$ can be expressed in terms of $P_{2}\times K_{5}$, which has a $(8; p, q)$-decomposition by Lemma \ref{LEM4C8P8}. Further $ K_{8,8}\times K_{3} $ and $K_{8,8}\times K_{5}$ has a $(8; p, q)$-decomposition, by Lemmas \ref{LEM0C8P8} and \ref{LEM3C8P8}. \\
	
	\noindent{\bf Subcase 1.6:} $n \equiv 1($mod $8)$.\\
	\noindent Let $K_{n}\cong \left(\frac{n-1}{8}\right) K_{9} \oplus \left(\frac{(n-1)(n-9)}{128}\right) K_{8,8}$ and $K_{m} \cong \left(\frac{m(m-1)}{2}\right) K_{2}$. Consider $K_{n} \times K_{m} \cong \left(\frac{m(m-1)(n-1)}{16}\right)  (K_{9} \times K_{2}) \oplus \left(\frac{m(m-1)(n-1)(n-9)}{256}\right)  (K_{8,8} \times K_{2})$. By Lemmas \ref{LEM7C8P8}, $K_{9} \times K_{2}$ has a $(8; p, q)$-decomposition. By Lemma \ref{LEM0C8P8} and \ref{LEM3C8P8}, $K_{8,8} \times K_{2}$ has a $(8; p, q)$-decomposition. \\
	
	\noindent Now $(8;p,q)$-decomposition of $K_{m} \times  K_{n}(\lambda)$ (where $\lambda$ is any non negative integers) follows for Case 1.\\
	
	\noindent{\bf Case 2:} $\lambda=2$.\\
	\noindent{\bf Subcase 2.1:} $n \equiv  0$(mod $4$) and $m \equiv  2$(mod $4$).\\
	Here $K_{n}\cong \left(\frac{n}{4}\right)  K_{4} \oplus \left(\frac{n^{2}-4n}{32}\right)  K_{4,4}$ and $K_{m}\cong K_{6} \oplus \left(\frac{m-6}{4}\right)  K_{4} \oplus \left(\frac{m-6}{4}\right)  K_{2,4} \oplus \left(\frac{(m-2)(m-6)}{32}\right)  K_{4,4}$. Consider $(K_{n} \times K_{m})(2) \cong \left(\frac{n}{4}\right)(K_{4} \times K_{6})(2) \oplus \left(\frac{n(m-6)}{16}\right) \\ (K_{4} \times K_{4})(2) \oplus \left(\frac{n(m-6)}{16}\right)  (K_{4} \times K_{2,4})(2) \oplus \left(\frac{n(m-6)(m+n-6)}{128}\right) (K_{4} \times K_{4,4})(2) \oplus \left(\frac{n(n-4)}{32}\right)  (K_{4,4} \times K_{6})(2) \oplus \left(\frac{n(n-4)(m-6)}{128}\right)  (K_{4,4} \times K_{2,4})(2) \oplus \left(\frac{n(n-4)(m-2)(m-6)}{1024}\right) \\(K_{4,4} \times K_{4,4})(2)$. By Lemma \ref{LEM1C8P8}, $K_{4} \times K_{4}(2)$ admits a $(8; p, q)$-decomposition. By Theorem \ref{THM5C8P8} and \ref{THM6C8P8}, $K_{4} \times K_{6}(2)$ and $K_{6} \times K_{4,4} (2)$ can be expressed in terms of $C_{4} \times K_{6}$ which has a $(8; p, q)$-decomposition by Lemma \ref{LEMC_{4}K_{6}C8P8}. The subgraphs $K_{4,4} \times K_{2,4}(2)$ and $K_{4,4} \times K_{4,4}(2)$ can be expressed in terms of $P_{2} \circ \overline{K}_{4}$ and hence admits a $(8; p, q)$-decomposition by Lemma \ref{LEM2C8P8}. By Theorem \ref{THM6C8P8}, $K_{4} \times K_{2,4}(2)$ and $K_{4} \times K_{4,4}(2)$ can be expressed in terms of $K_{4} \times C_{4}$ and hence admits a $(8; p, q)$-decomposition by Lemma \ref{LEMC_{4}K_{4}C8P8}. \\
	
	\noindent{\bf Subcase 2.2:} $n \equiv  1$(mod $4$) and $m \equiv  2$(mod $4$).\\
	Here $K_{m}\cong K_{6} \oplus \left(\frac{m-6}{4}\right)  K_{4} \oplus \left(\frac{m-6}{4}\right)  K_{2,4} \oplus \left(\frac{(m-2)(m-6)}{32}\right)  K_{4,4}$ and $K_{n}\cong \left(\frac{n-1}{4}\right)  K_{5} \oplus \left(\frac{(n-1)(n-5)}{32}\right)  K_{4,4}$. Consider {\footnotesize $(K_{n} \times K_{m})(2) \cong \left(\frac{n-1}{4}\right)  (K_{5} \times K_{6})(2) \oplus \left(\frac{(m-6)(n-1)}{16}\right)  (K_{5} \times K_{4})(2) \oplus \left(\frac{(m-6)(n-1)}{16}\right)  (K_{5} \times K_{2,4})(2) \oplus  \left(\frac{(m-2)(m-6)(n-1)}{128}\right) (K_{5} \times K_{4,4})(2) \oplus \left(\frac{(n-1)(n-5)}{32}\right)  (K_{4,4} \times K_{6})(2) \oplus  \left(\frac{(m-6)(n-1)(n-5)}{128}\right)  (K_{4,4} \times K_{4})(2) \oplus \left(\frac{(m-6)(n-1)(n-5)}{128}\right) \\ (K_{4,4} \times K_{2,4})(2) \oplus  \left(\frac{(m-2)(m-6)(n-1)(n-5)}{1024}\right)  (K_{4,4} \times K_{4,4})(2)$}. By Theorem \ref{THM6C8P8}, $K_{4,4} \times K_{4}(2)$ can be expressed in terms of $C_{4} \times K_{4}$ and hence admits a $(8; p, q)$-decomposition by Lemma \ref{LEMC_{4}K_{4}C8P8}. By Theorem \ref{THM5C8P8} and \ref{THM6C8P8}, $K_{5} \times K_{6}(2)$ and $K_{4,4} \times K_{6}(2)$ can be expressed in terms of $C_{4} \times K_{6}$ and hence admits a $(8; p, q)$-decomposition by Lemma \ref{LEMC_{4}K_{6}C8P8}. $K_{5} \times K_{4}(2), \ K_{5} \times K_{2,4}(2)$ and $K_{5} \times K_{4,4}(2)$ can be expressed in terms of $P_{2} \times K_{5}$ and hence admits a $(8; p, q)$-decomposition by Lemma \ref{LEM4C8P8}. The subgraphs $K_{4,4} \times K_{2,4}(2)$ and $K_{4,4} \times K_{4,4}(2)$ can be expressed in terms of $P_{2} \circ \overline{K}_{4}$ and hence admits a $(8; p, q)$-decomposition by Lemma \ref{LEM2C8P8}. \\
	
	\noindent{\bf Subcase 2.3:}	$n \equiv  0$(mod $4$) and $m \equiv  1$(mod $2$).\\
	Here $K_{n}\cong \left(\frac{n}{4}\right)  K_{4} \oplus \left(\frac{n^{2}-4n}{32}\right)  K_{4,4}$ and by Theorem \ref{THM2C8P8}, $K_{m}$ is decomposed into subgraphs isomorphic to $K_{3}$ and $K_{5}$. By Theorems \ref{THM5C8P8} and \ref{THM6C8P8}, $K_{4} \times K_{3}(2)$ and $K_{4,4} \times K_{3}(2)$ can be expressed in terms of $C_{4} \times K_{3}$ and hence admits a $(8; p, q)$-decomposition by Lemma \ref{LEMC_{4}K_{3}C8P8}. By Theorems \ref{THM1C8P8} and \ref{THM3C8P8}, $K_{4} \times K_{5}(2)$ and $K_{4,4} \times K_{5}(2)$ can be expressed in terms of $P_{2} \times K_{5}$ and hence admits a $(8; p, q)$-decomposition by Lemma \ref{LEM4C8P8}. \\

	\noindent{\bf Subcase 2.4:}	$n \equiv  1$(mod $4$) and $m \equiv  1$(mod $2$).\\
	Here $K_{n}\cong \left(\frac{n-1}{4}\right)  K_{5} \oplus \left(\frac{(n-1)(n-5)}{32}\right)  K_{4,4}$ and by Theorem \ref{THM2C8P8}, $K_{m}$ is decomposed into subgraphs isomorphic to $K_{3}$ and $K_{5}$. By Theorems \ref{THM5C8P8} and \ref{THM6C8P8}, $K_{5} \times K_{3}(2)$ and $K_{4,4} \times K_{3}(2)$ can be expressed in terms of $C_{4} \times K_{3}$ and hence admits a $(8; p, q)$-decomposition by Lemma \ref{LEMC_{4}K_{3}C8P8}. By Theorem \ref{THM1C8P8} and \ref{THM3C8P8}, $K_{5} \times K_{5}(2)$ and $K_{4,4} \times K_{5}(2)$ can be expressed in terms of $P_{2} \times K_{5}$ and hence admits a $(8; p, q)$-decomposition by Lemma \ref{LEM4C8P8}.\\
	
	\noindent{\bf Case 3:} $\lambda=4$.\\
	\noindent{\bf Subcase 3.1:} $n \equiv  1$(mod $2$) and $m \equiv  0$(mod $2$).\\
	\noindent{\bf Subcase 3.1.1:} $n \equiv  1$(mod $2$) and $m \equiv  0$(mod $4$).\\
	The proof of this case follows from Subcase 2.3.\\
	
	\noindent{\bf Subcase 3.1.2:} $n \equiv  1$(mod $2$) and $m \equiv  2$(mod $4$).\\
	Here $K_{m}\cong K_{6} \oplus \left(\frac{m-6}{4}\right)  K_{4} \oplus \left(\frac{m-6}{4}\right)  K_{2,4} \oplus \left(\frac{(m-2)(m-6)}{32}\right)  K_{4,4}$ and by Theorem \ref{THM2C8P8}, $K_{m}$ is decomposed into subgraphs isomorphic to $K_{3}$ and $K_{5}$. By Theorem \ref{THM1C8P8}, $K_{3} \times K_{6}(4)$ can be expressed in terms of $K_{3} \times P_{4}$ and hence admits a $(8; p, q)$-decomposition by Lemma \ref{LEM5C8P8}. By Theorems \ref{THM5C8P8} and \ref{THM6C8P8}, $K_{3} \times K_{4}(4),\ K_{3} \times K_{2,4}(4)$ and $K_{3} \times K_{4,4}(4)$ can be expressed in terms of $K_{3} \times C_{4}$ and hence admits a $(8; p, q)$-decomposition by Lemma \ref{LEMC_{4}K_{3}C8P8}. By Theorem \ref{THM1C8P8} and \ref{THM3C8P8}, $K_{5} \times K_{6}(4),\ K_{5} \times K_{4}(4),\ K_{5} \times K_{2,4}(4)$ and $K_{5} \times K_{4,4}(4)$ can be expressed in terms of $K_{5} \times P_{2}$ and hence admits a $(8; p, q)$-decomposition by Lemma \ref{LEM4C8P8}.  \\
	
	\noindent{\bf Subcase 3.2:} $n \equiv  0$(mod $2$) and $m \equiv  0$(mod $2$).\\
	\noindent{\bf Subcase 3.2.1:} $n \equiv  0$(mod $4$) and $m \equiv  0$(mod $4$).\\
	The proof of this case follows from Subcase 1.1.\\
	
	\noindent{\bf Subcase 3.2.2:} $n \equiv  2$(mod $4$) and $m \equiv  0$(mod $4$).\\
	The proof of this case follows from Subcase 2.1.\\
	
	\noindent{\bf Subcase 3.2.3:} $n \equiv  2$(mod $4$) and $m \equiv  2$(mod $4$).\\
	Here $K_{n}\cong K_{6} \oplus \left(\frac{n-6}{4}\right)  K_{4} \oplus \left(\frac{n-6}{4}\right)  K_{2,4} \oplus \left(\frac{(n-2)(n-6)}{32}\right)  K_{4,4}$ and $K_{m}\cong K_{6} \oplus \left(\frac{m-6}{4}\right)  K_{4} \oplus \left(\frac{m-6}{4}\right)  K_{2,4} \oplus \left(\frac{(m-2)(m-6)}{32}\right)  K_{4,4}$. Consider {\footnotesize $(K_{n} \times K_{m})(4) \cong (K_{6} \times K_{6})(4) \oplus \left(\frac{n+m-12}{4}\right)  (K_{6} \times K_{4})(4) \oplus \left(\frac{n+m-12}{4}\right) (K_{6} \times K_{2,4})(4) \oplus  \left(\frac{((n-2)(n-6))+((m-2)(m-6))}{32}\right) \\ (K_{6} \times K_{4,4})(4)  \oplus \left(\frac{(n-6)(m-6)(m+n-4)}{128}\right) (K_{4} \times K_{4,4})(4) \oplus \left(\frac{(m-6)(n-6)}{16}\right)(K_{4} \times K_{4})(4) \oplus \left(\frac{(m-2)(m-6)(n-2)(n-6)}{1024}\right)  (K_{4,4} \times K_{4,4})(4)  \oplus  \left(\frac{(m-6)(n-6)}{16}\right)(K_{2,4} \times K_{2,4})(4) \oplus \left(\frac{(m-6)(n-6)}{8}\right) \\ (K_{2,4} \times K_{4})(4) \oplus \left(\frac{(m-6)(n-6)(m+n-4)}{128}\right)  (K_{2,4} \times K_{4,4})(4) $}. By Lemma \ref{LEM1C8P8}, $K_{4} \times K_{4}(4)$ admits a $(8; p, q)$-decomposition. By Theorem \ref{THM5C8P8} and \ref{THM6C8P8}, $K_{6} \times K_{6}(4),\ K_{6} \times K_{4}(4),\ K_{6} \times K_{2,4}(4),\ K_{6} \times K_{4,4}(4)$ can be expressed in terms of $C_{4} \times K_{6}$ admits a $(8; p, q)$-decomposition by Lemma \ref{LEMC_{4}K_{6}C8P8}. By Theorem \ref{THM6C8P8}, $K_{4,4} \times K_{4}(4),\ K_{2,4} \times K_{4}(4)$ can be expressed in terms of $K_{4} \times C_{4}$ and hence admits a $(8; p, q)$-decomposition by Lemma \ref{LEMC_{4}K_{4}C8P8}. The subgraphs $K_{2,4} \times K_{2,4}(4),\ K_{2,4} \times K_{4,4}(4),\ K_{4,4} \times K_{4,4}(4)$ can be expressed in terms of $P_{2} \circ \overline{K}_{4}$ and hence admits a $(8; p, q)$-decomposition by Lemma \ref{LEM2C8P8}. \\
	
	\noindent{\bf Subcase 3.3:} $n \equiv  1$(mod $2$) and $m \equiv  1$(mod $2$).\\
	By Theorem \ref{THM2C8P8}, $K_{m}$ and $K_{n}$ is decomposed into subgraphs isomorphic to $K_{3}$ and $K_{5}$. By Lemma \ref{LEM4K_{3}K_{3}C8P8}, $K_{3} \times K_{3}(4)$ admits a $(8; p, q)$-decomposition. By Theorem \ref{THM1C8P8}, $K_{3} \times K_{5}(4),\ K_{5} \times K_{5}(4)$ can be expressed in terms of $P_{2} \times K_{5}$ admits a $(8; p, q)$-decomposition by Lemma \ref{LEM4C8P8}. \\
	\noindent This completes the proof of the theorem.\hfill$\Box$
	
	\section{$(8; p, q)$-decomposition of $K_{m} \circ \overline{K}_{n}(\lambda)$}
	In this section, we prove the necessary and sufficient conditions for $(8; p, q)$-decomposition of $ K_{m} \circ \overline{K}_{n}(\lambda)$.
	
	\begin{lemma}
		\label{LEM00C8P8}
		If $K_{m} \circ \overline{K}_{2}$ has a $(8; p, q)$-decomposition, then  $K_{m} \circ \overline{K}_{2t}$ has a $(8; p, q)$-decomposition.\\
	\end{lemma}
	\noindent{\bf Proof.}
	$K_{m} \circ \overline{K}_{2t} \cong (K_{m} \circ \overline{K}_{2})\circ \overline{K}_{t}\cong t (K_{m} \circ \overline{K}_{2}) \oplus \frac{t(t-1)}{2} ((K_{m} \circ \overline{K}_{2}) \times K_{2}) $. It is easy to see that $K_{m} \circ \overline{K}_{2t}$ has a $(8; p, q)$-decomposition when $K_{m} \circ \overline{K}_{2}$ has a $(8; p, q)$-decomposition.\hfill$\Box$

	\begin{lemma}
		Let $G$ be the graph such that $V(G)=A\cup B$, with $|A|=|B|=8$ and edges $E(G)= \langle A \rangle \cup \langle A,B \rangle$, where $\langle A \rangle \cong K_{8}-F_{1}$ and $\langle A,B \rangle \cong K_{8,8}$. Then there exists a $(8; p, q)$-decomposition of $G$.
		\label{LEM9C8P8}
	\end{lemma}
	\noindent{\bf Proof.}
	Here $G \cong K_{8,8}\oplus (K_{8}-F_{1})$. Let the vertices of the graph $G$ be partitioned into $A=\{a_{i}|0 \leq i \leq 7\}$ and $B=\{b_{i}|0 \leq i \leq 7\}$. The decomposition of subgraph $(K_{8}-F_{1})$ into $8$-cycles given by: $C^{1}=(a_{0}a_{7}a_{5}a_{4}a_{6}a_{3}a_{1}a_{2}),\ C^{2}=(a_{0}a_{1}a_{4}a_{7}a_{6}a_{5}a_{2}a_{3})$ and $C^{3}=(a_{0}a_{5}a_{3}a_{4}a_{2}a_{7}a_{1}a_{6})$. The decomposition of $ K_{8,8} $ into $8$-cycles given by: $C^{4}=(a_{0}b_{2}a_{6}b_{0}a_{4}b_{6}a_{2}b_{4}), \ C^{5}=(a_{0}b_{1}a_{2}b_{3}a_{4}b_{5}a_{6}\\ b_{7}),  C^{6}=(a_{0}b_{3}a_{6}b_{1}a_{4}b_{7}a_{2}b_{5})$ and $C^{7}=(a_{0}b_{0}a_{2}b_{2}a_{4}b_{4}a_{6}b_{6})$. Consider the permutation $\rho=(a_{0}a_{1}a_{2}a_{3}a_{4}a_{5}a_{6}a_{7})(b_{0}b_{1}b_{2}b_{3}b_{4}b_{5}b_{6}b_{7})$. This permutation acting on $C^{4},\ C^{5},\ C^{6}$ and $C^{7}$ provides remaining $8$-cycles of $K_{8,8}$, say $C^{8},\ C^{9},\ C^{10}$ and $C^{11}$ respectively. The required $(8; p, q)$-decomposition of $G$ is as follows:\\
	$D_{1}$: Consider $C^{4}\cup C^{5}$ and decompose them into two copies of $P_{8}$ as follows: $P^{1}= b_{1}a_{0}b_{2}a_{6}b_{0}a_{4}b_{6}a_{2}b_{4} $ and $ P^{2}= b_{1}a_{2}b_{3}a_{4}b_{5}a_{6}b_{7}a_{0}b_{4} $. \\
	$D_{2}$: Decomposition of $C^{8}\cup C^{9}$ into two copies of $P_{8}$ given by: $ \rho(P^{1}) $ and $\rho(P^{2}) $.\\
	$D_{3}$: Consider $C^{2}\cup C^{11}$ and decompose them into two copies of $P_{8}$ as follows: $P^{3}= a_{0}a_{3}a_{2}a_{5}a_{6}a_{7}a_{4}a_{1}b_{7} $ and $P^{4}=a_{0}a_{1}b_{1}a_{3}b_{3}a_{5}b_{5}a_{7}b_{7} $.\\
	$D_{4}$: Decompose $C^{3}\cup C^{10}$ and into two copies of $P_{8}$, say $P^{5}= a_{0}a_{5}b_{2}a_{7}b_{4}a_{1}b_{6}\\a_{3}b_{0} $ and $ P^{6}= a_{0}a_{6}a_{1}a_{7}a_{2}a_{4}a_{3}a_{5}b_{0} $. \\
	$D_{5}$: By decomposing $C^{1}\cup C^{7}$ we obtain two copies of $P_{8}$, say $ P^{7}= b_{0}a_{0}a_{2}a_{1}a_{3}\\a_{6}a_{4}a_{5}a_{7}$ and $ P^{8}=a_{7}a_{0}b_{6}a_{6}b_{4}a_{4}b_{2}a_{2}b_{0} $.\\
	$D_{6}$: Decomposition of $C^{1}\cup  C^{6}\cup C^{7}$ into three copies of $P_{8}$ given by: $ P^{9}=b_{0}a_{0}b_{3}a_{6}b_{1}a_{4}b_{7}a_{2}b_{5},\  P^{10}=b_{5}a_{0}a_{2}a_{1}a_{3}a_{6}a_{4}a_{5}a_{7}$ and $P^{8}=a_{7}a_{0}b_{6}a_{6}b_{4}a_{4}b_{2}a_{2}b_{0} $. \\
	Here $ (p,q) \in \{(0,11),(2,9),(3,8),(4,7),(5,6),(6,5),(7,4),(8,3),(9,2),(10,1),\\(11,0)\}$. For $1 \leq i \leq 6$, combinations of $D_{i}$ provide the required $(8; p, q)$-decomposition.\hfill$\Box$
	
	\begin{lemma}
		For $k\geq 2, K_{4k,4k}$ has a $(8; p, q)$-decomposition.
		\label{LEM11C8P8}
	\end{lemma}
	\noindent{\bf Proof.}
	For any $k\geq 2, K_{4k,4k}$ can be expressed as sum of $K_{8,8},\ K_{8,12}$ and $K_{12,12}$. In Lemma \ref{LEM3C8P8}, a $(8; p, q)$-decomposition of $K_{8,8}$ is discussed. Let $A=\{a_{i}|0 \leq i \leq 11\}$ and $B=\{b_{i}|0 \leq i \leq 11\}$. We partition the vertices as follows: $A_{1}=\{ a_{0}, a_{1}, a_{2}, a_{3}\},\ A_{2}=\{ a_{4}, a_{5}, a_{6}, a_{7}\},\ A_{3}=\{ a_{8}, a_{9}, a_{10}, a_{11}\},\ B_{1}=\{ b_{0}, b_{1}, b_{2}, b_{3}\},\ B_{2}=\{ b_{4}, b_{5}, b_{6}, b_{7}\}$ and $B_{3}=\{ b_{8}, b_{9},\\ b_{10}, b_{11}\}$. The subgraph induced by $\langle A_{i}, B_{i}\rangle \cong K_{4,4}, i=\{1,2,3\}$. Now, we first decompose $\{\langle A_{1}, B_{1}\rangle \cup \langle A_{1}, B_{2}\rangle \cup \langle A_{1}, B_{3}\rangle\}$ as follows:\\
	The decomposition of $\langle A_{1}, B_{1}\rangle$ into $8$-cycles given by: $C^{1}=(a_{0}b_{0}a_{3}b_{3}a_{2}b_{2}a_{1}b_{1})$ and $C^{2}=(a_{0}b_{2}a_{3}b_{1}a_{2}b_{0}a_{1}b_{3})$. Replacing $b_{0},\ b_{1},\ b_{2}$ and $b_{3}$ with $b_{4},\ b_{5},\ b_{6}$ and $b_{7}$ respectively $(b_{8},\ b_{9},\ b_{10}$ and $b_{11}$ respectively) we obtain a $C_{8}$-decomposition of $\langle A_{1}, B_{2}\rangle(\langle A_{1}, B_{3}\rangle)$, say $C^{3}$ and $C^{4}$ respectively($C^{5}$ and $C^{6}$respectively). \\
	$D_{1}$: Decomposition of $C^{1} \cup C^{4}$ into $P_{8}$ given by: $b_{0}a_{0}b_{6}a_{3}b_{5}a_{2}b_{4}a_{1}b_{7}$ and $b_{7}a_{0}b_{1}a_{1}b_{2}a_{2}b_{3}a_{3}b_{0}$.\\ $D_{2}$: Similarly we obtain a $P_{8}$-decomposition of $ C^{2} \cup C^{5}$ and $ C^{3} \cup C^{6} $, since $  C^{2} \cup C^{5} $ and $ C^{3} \cup C^{6} $ is isomorphic to $ C^{1} \cup C^{4} $. \\
	$D_{3}$: Consider $C^{1}\cup C^{2}\cup C^{4}$ and decompose them into three copies of $P_{8}$ as follows: $b_{7}a_{1}b_{2}a_{2}b_{3}a_{3}b_{0}a_{0}b_{1}$,\ $b_{1}a_{1}b_{4}a_{2}b_{5}a_{3}b_{6}a_{0}b_{2}$ and $b_{2}a_{3}b_{1}a_{2}b_{0}a_{1}b_{3}a_{0}b_{7}$. \\
	For $1 \leq i \leq 3$, combinations of $D_{i}$ provide the required $(8; p, q)$-decomposition of $\{\langle A_{1}, B_{1}\rangle \cup \langle A_{1}, B_{2}\rangle \cup \langle A_{1}, B_{3}\rangle\}$. Similarly we decompose $\{\langle A_{2}, B_{1}\rangle \cup \langle A_{2}, B_{2}\rangle \cup \langle A_{2}, B_{3}\rangle\}$ and $\{\langle A_{3}, B_{1}\rangle \cup \langle A_{3}, B_{2}\rangle \cup \langle A_{3}, B_{3}\rangle\}$. The decomposition of $\{\langle A_{i}, B_{j}\rangle, i=1,2$ and $j=1,2,3\}$ provides a $(8; p, q)$-decomposition of $K_{8,12}$ and the decomposition of $\{\langle A_{i}, B_{j}\rangle, i=1,2,3$ and $j=1,2,3\}$ provides a $(8; p, q)$-decomposition of $K_{12,12}$. \hfill$\Box$
	
	\begin{lemma}
		There exists a $(8; p, q)$-decomposition of $K_{10}-F_{1}$.
		\label{LEM12C8P8}
	\end{lemma}
	\noindent{\bf Proof.}
	Let the vertices of the graph $G$ be partitioned into $A=\{a_{i}|0 \leq i \leq 4\}$ and $B=\{b_{i}|0 \leq i \leq 4\}$. The decomposition of $K_{10}-F_{1}$ into $8$-cycles given by: $C^{1}=(a_{0}b_{1}a_{3}b_{4}a_{2}b_{3}a_{1}b_{2}),\ C^{2}=(a_{0}a_{2}a_{4}b_{2}b_{1}b_{4}b_{0}a_{3}),\ C^{3}=(a_{0}a_{1}a_{3}b_{2}b_{4}b_{3}b_{0}a_{4}),\ C^{4}=(a_{1}a_{4}a_{3}a_{2}b_{1}b_{3}b_{2}b_{0})$ and $C^{5}=(a_{0}b_{3}a_{4}b_{1}b_{0}a_{2}a_{1}b_{4})$. The required $(8; p, q)$-decomposition of $K_{10}-F_{1}$ is as follows:\\
	$D_{1}$: Consider $C^{1}\cup C^{2}$ and decompose them into $P_{8}$ as follows: $ P^{1}=b_{0}b_{4}b_{1}a_{3}a_{0}b_{2}a_{4}a_{2}b_{3} $ and $ P^{2}=b_{0}a_{3}b_{4}a_{2}a_{0}b_{1}b_{2}a_{1}b_{3} $. \\
	$D_{2}$: Decomposition $C^{4}\cup C^{5}$ into $P_{8}$ given by: $ P^{3}=b_{4}a_{0}b_{3}a_{4}a_{1}a_{2}b_{1}b_{0}b_{2} $ and $ P^{4}=b_{4}a_{1}b_{0}a_{2}a_{3}a_{4}b_{1}b_{3}b_{2} $.\\
	$D_{3}$: Consider $C^{1}\cup C^{2}\cup C^{3}$ and decompose them into $P_{8}$, say $ P^{1}=b_{0}b_{4}b_{1}a_{3}a_{0}\\b_{2} a_{4}a_{2}b_{3},\ P^{5}=b_{0}a_{3}b_{4}b_{3}a_{1}b_{2}b_{1}a_{0}a_{2} $ and $ P^{6}=a_{2}b_{4}b_{2}a_{3}a_{1}a_{0}a_{4}b_{0}b_{3} $.\\
	Here $ (p,q) \in \{(0,5),(2,3),(3,2),(4,1),(5,0)\}$. For $1 \leq i \leq 3$, combinations of $D_{i}$ provide the required $(8; p, q)$-decomposition.\hfill$\Box$
	
	\begin{lemma}
		The graph $(K_{12}-F_{1})(2)$ admits a $(8; p, q)$-decomposition.
		\label{LEM2K_{12}-1C8P8} 
	\end{lemma}
	\noindent{\bf Proof.}
	Let the vertices of the graph $(K_{12}-F_{1})(2)$ be partitioned into $A=\{a_{i}|0 \leq i \leq 5\}$ and $B=\{b_{i}|0 \leq i \leq 5\}$. The decomposition of $(K_{12}-F_{1})(2)$ into $8$-cycles given by: $C^{1}=(a_{0}a_{1}b_{2}a_{5}a_{3}a_{4}b_{4}a_{2}),\ C^{2}=(a_{1}b_{2}b_{0}b_{1}a_{5}b_{5}b_{3}b_{4}),\ C^{3}=(a_{0}a_{1}b_{2}a_{5}a_{4}a_{2}b_{1}a_{3})$ and the permutation $\rho=(a_{0}a_{1}a_{2}a_{3}a_{4}a_{5}a_{6}a_{7}a_{8})(b_{0}b_{1}b_{2}b_{3}b_{4}b_{5}b_{6}b_{7}b_{8})$. Then $\{\rho^{i} (C^{j}) | 0\leq i \leq 4 ,\ j=1,2,3\}$ provides a $C_{8}$-decomposition of $(K_{12}-F_{1})(2)$. The required $(8; p, q)$-decomposition of $(K_{12}-F_{1})(2)$ is as follows:\\
	$D_{1}$: Decomposition of $C^{1}\cup C^{2}$ into $2$ copies of $ P_{8}$ given by: $P^{1}=b_{1}  a_{5} b_{2}  a_{1}  a_{0}  a_{2} \\ b_{4} b_{3} b_{5}$ and $P^{2}=b_{5}  a_{5} a_{3}  a_{4}  b_{4}  a_{1}  b_{2} b_{0} b_{1}$.\\
	$D_{2}$: For $1\leq i \leq 4$, $\rho^{i}(P^{1})$ and $\rho^{i}(P^{2})$ provides a $P_{8}$-decomposition of $\rho^{i}(C^{1})\cup \rho^{i}(C^{2})$.\\
	$D_{3}$: Decompose $C^{1}\cup C^{2}\cup C^{3}$ into $P_{8}$ as follows: $P^{1}=b_{1}  a_{5} b_{2}  a_{1}  a_{0}  a_{2}  b_{4} b_{3} b_{5},\\ P^{3}=b_{5}  a_{5} b_{2}  a_{1}  b_{4}  a_{4}  a_{3} b_{1} a_{2}$ and $P^{4}=a_{2}  a_{4} a_{5}  a_{3}  a_{0}  a_{1}  b_{2} b_{0} b_{1}$.\\
	$D_{4}$: For $1\leq i \leq 4$, $\rho^{i}(P^{1}),\ \rho^{i}(P^{2})$ and $\rho^{i}(P^{3})$ provides a $P_{8}$-decomposition of $\rho^{i}(C^{1})\cup \rho^{i}(C^{2}) \cup \rho^{i}(C^{3})$.\\
	Here $(p,q) \in \{(0,15),(2,13),(3,12),(4,11),(5,10),(6,9),(7,8),(8,7),(9,6),\\(10,5),(11,4),(12,3),(13,2),(14,1),(15,0)\}$. For $1 \leq i \leq 4$, combinations of $D_{i}$ provide the required $(8; p, q)$-decomposition.\hfill$\Box$
	
	\begin{lemma}
		There exists a $(8; p, q)$-decomposition of $(K_{16}-F_{1})$.
		\label{LEM8C8P8}
	\end{lemma}
	\noindent{\bf Proof.}
	Consider $(K_{16}-F_{1}) \cong 2 K_{8}\oplus (K_{8,8}-F_{1})$. Partition the vertices of the graph into $A=\{a_{i}|0 \leq i \leq 7\}$ and $B=\{b_{i}|0 \leq i \leq 7\}$ such that $\langle A \rangle \cong \langle B \rangle \cong K_{8}$ and $ \langle A,B \rangle \cong (K_{8,8}-F_{1})$. Let $C^{1}=(a_{0}a_{3}a_{6}a_{1}a_{4}a_{7}a_{2}a_{5}), C^{2}=(a_{0}a_{6}a_{5}a_{1}a_{2}a_{4}a_{3}a_{7})$ and $C^{3}=(a_{0}a_{1}a_{7}a_{6}a_{2}a_{3}a_{5}a_{4})$. Replace $a_{i}$ by $b_{i}$ in the cycles $c_{j},\ 1 \leq j \leq 3$ and name them as $C^{4},\ C^{5}$ and $C^{6}$ respectively. Using the edges of $\langle A,B \rangle$, let us consider the $8$-cycles: $C^{7}=(a_{0}b_{1}a_{6}b_{7}a_{4}b_{5}a_{2}b_{3}),\\ C^{8}=(a_{0}b_{2}a_{6}b_{0}a_{4}b_{6}a_{2}b_{4})$ and $C^{9}=(a_{0}b_{5}a_{6}b_{3}a_{4}b_{1}a_{2}b_{7})$. Consider the permutation $\rho=(a_{0}a_{1}a_{2}a_{3}a_{4}a_{5}a_{6}a_{7})(b_{0}b_{1}b_{2}b_{3}b_{4}b_{5}b_{6}b_{7})$ acting on $C^{7},\ C^{8}$ and $C^{9}$ provides $8$-cycles, say $C^{10},\ C^{11}$ and $C^{12}$, respectively. The remaining edges form two $8$-cycles given by: $C^{13}=(a_{1}b_{7}b_{5}a_{7}a_{5}b_{3}b_{1}a_{3})$ and $C^{14}=(a_{0}b_{6}b_{4}a_{6}a_{4}b_{2}b_{0}a_{2})$. The required $(8; p, q)$-decomposition of $(K_{16}-F_{1})$ is as follows:\\
	$D_{1}$: Consider $C^{1}\cup C^{13}$ and decompose them into two copies of $P_{8}$ as follows: $P^{1}= a_{0}a_{3}a_{6}a_{1}a_{4}a_{7}a_{2}a_{5}b_{3} $ and $ P^{2}=a_{0}a_{5}a_{7}b_{5}b_{7}a_{1}a_{3}b_{1}b_{3} $.\\
	$D_{2}$: By decomposing $C^{2}\cup C^{14}$ we get two copies of $P_{8}$, say $P^{3}= b_{4}b_{6}a_{0}a_{2}b_{0}b_{2}\\a_{4}a_{6}a_{5} $ and $ P^{4}=b_{4}a_{6}a_{0}a_{7}a_{3}a_{4}a_{2}a_{1}a_{5} $.\\
	$D_{3}$: Decomposition of $C^{3}\cup C^{7}$ into two copies of $P_{8}$ given by: $ P^{5}=b_{5}a_{4}a_{0}a_{1}a_{7}a_{6}a_{2}a_{3}a_{5} $ and $ P^{6}=a_{5}a_{4}b_{7}a_{6}b_{1}a_{0}b_{3}a_{2}b_{5} $.\\
	$D_{4}$: Consider $C^{4}\cup C^{8}$ and decompose them into two copies of $P_{8}$ as follows: $ P^{7}=a_{6}b_{0}b_{5}b_{2}b_{7}b_{4}b_{1}b_{6}b_{3} $ and $ P^{8}=a_{6}b_{2}a_{0}b_{4}a_{2}b_{6}a_{4}b_{0}b_{3} $.\\
	$D_{5}$: Decompose $C^{5}\cup C^{9}$ into two copies of $P_{8}$, say $ P^{9}=b_{2}b_{4}b_{3}b_{7}b_{0}b_{6}b_{5}b_{1}a_{2} $ and $ P^{10}=a_{2}b_{7}a_{0}b_{5}a_{6}b_{3}a_{4}b_{1}b_{2} $.\\
	$D_{6}$: Decomposition of $C^{6}\cup C^{10}$ into two copies of $P_{8}$ given by: $ P^{11}=a_{1}b_{4}b_{0}b_{1}b_{7}b_{6}b_{2}b_{3}b_{5} $ and $ P^{12}=b_{5}b_{4}a_{3}b_{6}a_{5}b_{0}a_{7}b_{2}a_{1} $.\\
	$D_{7}$: By decomposing $C^{11}\cup C^{12}$ we obtain two copies of $P_{8}$, say $P^{13}= b_{3}a_{1}b_{0}a_{3}b_{2}a_{5}b_{4}a_{7}b_{6} $ and $ P^{14}=b_{3}a_{7}b_{1}a_{5}b_{7}a_{3}b_{5}a_{1}b_{6} $.\\
	$D_{8}$: Decompose $C^{1}\cup C^{7}\cup C^{13}$ into three copies of $P_{8}$ as follows: $P^{2}= a_{0}a_{5}a_{7}b_{5}b_{7}a_{1}a_{3}b_{1}b_{3},\ P^{15}=a_{0}a_{3}a_{6}a_{1}a_{4}a_{7}a_{2}b_{3}a_{5} $ and $ P^{16}=a_{5}a_{2}b_{5}a_{4}b_{7}a_{6}b_{1}a_{0}b_{3} $.\\
	Here $(p,q) \in \{(0,14),(2,12),(3,11),(4,10),(5,9),(6,8),(7,7),(8,6),(9,5),\\(10,4),(11,3),(12,2)(13,1),(14,0)\}$. For $1 \leq i \leq 8$, combinations of $D_{i}$ provide the required $(8; p, q)$-decomposition.\hfill$\Box$
	
	\begin{lemma}
		The graph $K_{4,6}$ admits a $(8; p, q)$-decomposition.
		\label{LEMK_{4,6}C8P8}
	\end{lemma}
	\noindent{\bf Proof.}
	Let the vertices of the graph $K_{4,6}$ be partitioned into $A=\{a_{i}|0 \leq i \leq 3\}$ and $B=\{b_{i}|0 \leq i \leq 5\}$. The decomposition of $K_{4,6}$ into $8$-cycles given by: $C^{1}=(a_{0}b_{2}a_{2}b_{4}a_{1}b_{3}a_{3}b_{0}),\ C^{2}=(a_{0}b_{1}a_{2}b_{3}a_{3}b_{5}a_{1}b_{4})$ and $C^{3}=(a_{0}b_{1}a_{2}b_{0}a_{1}b_{5}a_{3}b_{2})$. The required $(8; p, q)$-decomposition of $K_{4,6}$ is as follows:\\
	For $(p,q)=(2,1)$: Consider $C^{1}\cup C^{2}$ and decompose them into $2$ copies of $ P_{8}$ as follows: $P^{1}=b_{1} a_{0} b_{2} a_{2} b_{4} a_{1} b_{3} a_{3} b_{0}$ and $P^{2}=b_{0} a_{0} b_{4} a_{1} b_{5} a_{3} b_{3} a_{2} b_{1}$.\\
	For $(p,q)=(3,0)$: Decomposition of $C^{1}\cup C^{2}\cup C^{3}$ into $3$ copies of $ P_{8}$ given by: $P^{1}=b_{1} a_{0} b_{2} a_{2} b_{4} a_{1} b_{3} a_{3} b_{0},\ P^{3}=b_{1} a_{2} b_{0} a_{1} b_{5} a_{3} b_{2} a_{0} b_{4}$ and $P^{4}=b_{4} a_{1} b_{5} a_{3} b_{3} a_{2} b_{1} a_{0} b_{0}$. \hfill$\Box$
	
	\begin{lemma}
		The graph $K_{6,6}(2)$ admits a $(8; p, q)$-decomposition.
		\label{LEM2K_{6,6}C8P8}
	\end{lemma}
	\noindent{\bf Proof.}
	Let the vertices of the graph $K_{6,6}(2)$ be partioned into $A=\{a_{i}|0 \leq i \leq 5\}$ and $B=\{b_{i}|0 \leq i \leq 5\}$. The decomposition of $K_{6,6}(2)$ into $8$-cycles given by: $C^{1}=(a_{0}b_{2}a_{2}b_{4}a_{4}b_{5}a_{5}b_{0}),\ C^{1'}=(a_{0}b_{2}a_{2}b_{4}a_{4}b_{5}a_{5}b_{0}),\ C^{2}=(a_{1}b_{1}a_{4}b_{0}a_{2}b_{5}a_{3}b_{3}),\
	C^{3}=(a_{0}b_{3}a_{5}b_{2}a_{3}b_{5}a_{1}b_{4}),\
	C^{4}=(a_{0}b_{1}a_{5}b_{4}a_{3}b_{0}a_{1}b_{5}),\ C^{5}=(a_{1}b_{0}a_{3}b_{1}a_{2}b_{3}a_{4}b_{2}),\
	C^{6}=(a_{0}b_{5}a_{2}b_{3}a_{1}b_{1}a_{5}b_{4}),\
	C^{7}=(a_{2}b_{0}a_{4}b_{3}a_{5}b_{2}a_{3}b_{1})$ and $C^{8}=(a_{0}b_{1}a_{4}b_{2}a_{1}b_{4}a_{3}b_{3})$. The required $(8; p, q)$-decomposition of $K_{6,6}(2)$ is as follows:\\
	$D_{1}$: Consider $C^{1}\cup C^{2}$ and decompose them into $2$ copies of $ P_{8}$ as follows: $P^{1}=b_{1} a_{4} b_{5} a_{5} b_{0} a_{0} b_{2} a_{2} b_{4}$ and $P^{2}=b_{4} a_{4} b_{0} a_{2} b_{5} a_{3} b_{3} a_{1} b_{1}$.\\
	$D_{2}$: Decomposition of $C^{3}\cup C^{4}$ into $2$ copies of $ P_{8}$ given by: $P^{3}=b_{1} a_{0} b_{4} a_{1} b_{5} a_{3} \\ b_{2} a_{5} b_{3}$ and $P^{4}=b_{3} a_{0} b_{5} a_{1} b_{0} a_{3} b_{4} a_{5} b_{1}$.\\
	$D_{3}$: Decomposition of $C^{5}\cup C^{6}$ into $2$ copies of $ P_{8}$ given by: $P^{5}= a_{5} b_{1} a_{2} b_{3} a_{4} b_{2} \\ a_{1} b_{0} a_{3}$ and $P^{6}=a_{3} b_{1} a_{1} b_{3} a_{2} b_{5} a_{0} b_{4} a_{5} $. \\
	$D_{4}$: Decomposition of $C^{7}\cup C^{8}$ into $2$ copies of $ P_{8}$ given by: $P^{7}= b_{0} a_{2} b_{1} a_{3} b_{2} a_{5} \\ b_{3} a_{4} b_{1}$ and $P^{8}= b_{1} a_{0} b_{3} a_{3} b_{4} a_{1} b_{2} a_{4} b_{0}$.\\
	$D_{5}$: Decomposition of $C^{1'}\cup C^{7}\cup C^{8}$ into $3$ copies of $ P_{8}$ given by: $P^{7}=b_{0} a_{2} b_{1} a_{3} b_{2} a_{5} b_{3} a_{4} b_{1}, P^{9}=b_{1} a_{0} b_{0} a_{4} b_{2} a_{2} b_{4} a_{3} b_{3}$ and $ P^{10}=b_{3} a_{0} b_{2} a_{1} b_{4} a_{4} b_{5} a_{5} b_{0}$.\\
	Here $(p,q) \in \{(0,9),(2,7),(3,6),(4,5),(5,4),(6,3),(7,2),(8,1),(9,0)\}$. For $1 \leq i \leq 5$, combinations of $D_{i}$ provide the required $(8; p, q)$-decomposition.\hfill$\Box$
	
	\begin{lemma}
		The graph $K_{4,10}$ admits a $(8; p, q)$-decomposition.
		\label{LEMK_{4,10}C8P8}
	\end{lemma}
	\noindent{\bf Proof.}
	Partition the vertices of the graph $K_{4,10}$ into $A=\{a_{i}|0 \leq i \leq 3\}$ and $B=\{b_{i}|0 \leq i \leq 9\}$. It is easy to see that, $K_{4,10} \cong K_{4,6} \oplus K_{4,4}$.  The subgraph $K_{4,6}$ is decomposed into $3 C_{8}$, say $C^{1},\ C^{2}$ and $C^{3}$ which is given in Lemma \ref{LEMK_{4,6}C8P8}. The decomposition of $K_{4,4}$ into $8$-cycles are given by: $C^{4}=(b_{6}a_{0}b_{7}a_{1}b_{8}a_{2}b_{9}a_{3})$ and $C^{5}=(a_{0}b_{8}a_{3}b_{7}a_{2}b_{6}a_{1}b_{9})$. The required $(8; p, q)$-decomposition of $K_{4,10}$ is as follows:\\
	$D_{1}$: Consider $C^{1}\cup C^{4}$ and decompose them into $2$ copies of $ P_{8}$ as follows: $P^{1}=b_{8} a_{1} b_{3} a_{3} b_{0} a_{0} b_{2} a_{2} b_{4}$ and $P^{2}=b_{4} a_{1} b_{7} a_{0} b_{6} a_{3} b_{9} a_{2} b_{8}$.\\
	$D_{2}$: Decomposition of $C^{2}\cup C^{5}$ into $2$ copies of $ P_{8}$ given by: $P^{3}=b_{9} a_{1} b_{4} a_{0} b_{1} a_{2} \\ b_{3} a_{3} b_{5}$ and $P^{4}=b_{5} a_{1} b_{6} a_{2} b_{7} a_{3} b_{8} a_{0} b_{9}$.\\
	$D_{3}$: Decomposition of $C^{2}\cup C^{3}\cup C^{5}$ into $3$ copies of $ P_{8}$ given by: $P^{3}=b_{9} a_{1} b_{4} a_{0} b_{1} a_{2} b_{3} a_{3} b_{5},\ P^{5}=b_{9} a_{0} b_{2} a_{3} b_{5} a_{1} b_{0} a_{2} b_{1}$ and $P^{6}=b_{1} a_{0} b_{8} a_{3} b_{7} a_{2} b_{6} a_{1} b_{5}$. \\
	Here $ (p,q) \in \{(0,5),(2,3),(3,2),(4,1),(5,0)\}$. For $1 \leq i \leq 3$, combinations of $D_{i}$ provide the required $(8; p, q)$-decomposition.\hfill$\Box$
	
	\begin{lemma}
		There exists a $(8; p, q)$-decomposition of $K_{18}-F_{1}$.
		\label{LEMK_{18}-1C8P8}
	\end{lemma}
	\noindent{\bf Proof.}
	$K_{18}- \textrm{F}_{1} \cong K_{10}- \textrm{F}_{1} \oplus K_{8}- \textrm{F}_{1} \oplus K_{4,4} \oplus K_{4,6} \oplus K_{4,10}$. Let the vertices of the graph $K_{18}-F_{1}$ be partitioned into $A=\{a_{i}|0 \leq i \leq 7\}$ and $B=\{b_{i}|0 \leq i \leq 9\}$. The subgraphs $K_{10}- \textrm{F}_{1},K_{4,6}$ and $K_{4,10}$ admits a $(8; p, q)$-decomposition by Lemmas \ref {LEM12C8P8}, \ref{LEMK_{4,6}C8P8} and \ref{LEMK_{4,10}C8P8}. The cycles $C^{1},\ C^{2}$ and $C^{3}$  in Lemma \ref{LEM9C8P8} gives $8$-cycles decomposition of $K_{8}-F_{1}$ and the cycles $C^{4}$ and $C^{5}$ in Lemma \ref{LEMK_{4,10}C8P8} gives $8$-cycles decomposition of $K_{4,4}$. The required $(8; p, q)$-decomposition of $K_{4,10}$ is as follows:\\
	$D_{1}$: Consider $C^{1}\cup C^{4}$ and decompose them into $2$ copies of $ P_{8}$ as follows: $P^{1}=b_{0} a_{0} a_{2} a_{1} a_{3} a_{6} a_{4} a_{5} a_{7}$ and $P^{2}=a_{7} a_{0} b_{1} a_{1} b_{2} a_{2} b_{3} a_{3} b_{0}$.\\
	$D_{2}$: Decomposition of $C^{2}\cup C^{5}$ into $2$ copies of $ P_{8}$ given by: $P^{3}=b_{0} a_{1} a_{0} a_{3} a_{2} a_{5} \\ a_{6} a_{7} a_{4}$ and $P^{4}=a_{4} a_{1} b_{3} a_{0} b_{2} a_{3} b_{1} a_{2} b_{0}$.\\
	$D_{3}$: Decomposition of $C^{1}\cup C^{3}\cup C^{4}$ into $3$ copies of $ P_{8}$ given by: $P^{1}=b_{0} a_{0} a_{2} a_{1} a_{3} a_{6} a_{4} a_{5} a_{7},\ P^{5}=b_{0} a_{3} a_{4} a_{2} a_{7} a_{1} a_{6} a_{0} a_{5}$ and $P^{6}=a_{5} a_{3} b_{3} a_{2} b_{2} a_{1} b_{1} a_{0} a_{7}$. \\
	Here $ (p,q) \in \{(0,5),(2,3),(3,2),(4,1),(5,0)\}$. For $1 \leq i \leq 3$, combinations of $D_{i}$ provide the required $(8; p, q)$-decomposition.\hfill$\Box$
	
	\begin{lemma}
			For any even integer $n \geq 2$, $(K_{m}\circ \overline{K}_{n})(2)$ admits a $(8; p, q)$-decomposition.
		\label{LEM2K_{m}{K}_{2}C8P8}
	\end{lemma}
	\noindent{\bf Proof.}
		Since $n$ is even, $(K_{m}\circ \overline{K}_{n})(2)$ can be expressed in terms of $(K_{m}\circ \overline{K}_{2})(2)$ by Lemma \ref{LEM00C8P8}. It is easy to see that $(K_{m}\circ \overline{K}_{2})(2)$ is isomorphic to $(K_{2m}-F_{1})(2)$.\\
		
	\noindent When $ m \equiv 0($mod $8), (K_{2m}- \textrm{F}_{1})(2) \cong (K_{16k}- \textrm{F}_{1})(2) $ which will be partitioned into union of subgraphs isomorphic to $(K_{16}- \textrm{F}_{1})(2),\ K_{4,6}(2),\ K_{6,6}(2)$ and $K_{4,10}(2)$. Hence by Lemmas \ref{LEM8C8P8}, \ref{LEMK_{4,6}C8P8}, \ref{LEM2K_{6,6}C8P8} and \ref{LEMK_{4,10}C8P8}, $(K_{16k}- \textrm{F}_{1})(2)$ admits a $(8; p, q)$-decomposition.\\
	
	\noindent When $ m \equiv 1($mod $8), (K_{2m}- \textrm{F}_{1})(2) = (K_{16k+2}- \textrm{F}_{1}) (2)$ which will be partitioned into union of subgraphs isomorphic to $(K_{18}- \textrm{F}_{1})(2),\ (K_{16}- \textrm{F}_{1})(2),\ K_{4,6}(2),\ K_{6,6}(2)$ and $K_{4,10}(2)$. Hence by Lemmas \ref{LEM8C8P8}, \ref{LEMK_{4,6}C8P8}, \ref{LEM2K_{6,6}C8P8}, \ref{LEMK_{4,10}C8P8} and \ref{LEMK_{18}-1C8P8}, $(K_{16k+2}- \textrm{F}_{1})(2)$ admits a $(8; p, q)$-decomposition.\\
	
	\noindent When $ m \equiv 3($mod $8), (K_{2m}- \textrm{F}_{1})(2) = (K_{16k+6}- \textrm{F}_{1})(2) $ which will be partitioned into union of subgraphs isomorphic to ($K_{10}- \textrm{F}_{1})(2),\ ( K_{12}- \textrm{F}_{1})(2),\ (K_{16}- \textrm{F}_{1})(2),\ K_{4,6}(2),\ K_{6,6}(2)$ and $K_{4,10}(2)$. Hence by Lemmas \ref{LEM12C8P8}, \ref{LEM2K_{12}-1C8P8}, \ref{LEM8C8P8}, \ref{LEMK_{4,6}C8P8}, \ref{LEM2K_{6,6}C8P8} and \ref{LEMK_{4,10}C8P8}, $(K_{16k+6}- \textrm{F}_{1})(2)$ admits a $(8; p, q)$-decomposition.\\
	
	\noindent When $ m \equiv r($mod $8)$, where $r \in \{2,5,7\}.\ ( K_{2m}- \textrm{F}_{1})(2)=(K_{16k+2r}-\textrm{F}_{1})(2) $ which will be partitioned into union of subgraphs isomorphic to $(K_{10}- \textrm{F}_{1})(2),\ (K_{16}- \textrm{F}_{1})(2),\ K_{4,6}(2),\ K_{6,6}(2)$ and $K_{4,10}(2)$. Hence by Lemmas \ref{LEM12C8P8}, \ref{LEM8C8P8}, \ref{LEMK_{4,6}C8P8}, \ref{LEM2K_{6,6}C8P8} and \ref{LEMK_{4,10}C8P8}, $(K_{16k+2r}-\textrm{F}_{1})(2)$ admits a $(8; p, q)$-decomposition.\\
	
	\noindent When $ m \equiv r($mod $8),$ where $r \in \{4,6\}.\ (K_{2m}- \textrm{F}_{1})(2)=(K_{16k+2r}-\textrm{F}_{1})(2)$ which will be partitioned into union of subgraphs isomorphic to $(K_{12}- \textrm{F}_{1})(2),\ (K_{16}- \textrm{F}_{1})(2),\ K_{4,6}(2),\ K_{6,6}(2)$ and $K_{4,10}(2)$. Hence by Lemmas \ref{LEM2K_{12}-1C8P8}, \ref{LEM8C8P8}, \ref{LEMK_{4,6}C8P8}, \ref{LEM2K_{6,6}C8P8} and \ref{LEMK_{4,10}C8P8}, $(K_{16k+2r}- \textrm{F}_{1})(2)$ admits a $(8; p, q)$-decomposition.\hfill$\Box$
	
	\begin{lemma}
		For any even integers $n \geq 2$, $C_{8} \circ \overline{K}_{n}$ admits a $(8; p, q)$-decomposition.
		\label{LEMC_{8}K_{2}C8P8}
	\end{lemma}
	\noindent{\bf Proof.}
	By Lemma \ref{LEM00C8P8}, it is enough to consider $C_{8} \circ \overline{K}_{2}$. Let the vertices of $C_{8} \circ \overline{K}_{2}$ be partitioned into $A=\{a_{i}| 0 \leq i \leq 7\}$ and $B=\{b_{i}| 0 \leq i \leq 7\}$ such that $\langle A \rangle \cong \langle B \rangle\cong C_{8}$ and $\langle A,B \rangle \cong C_{8} \times K_{2}$.
	The decomposition of $C_{8} \circ \overline{K}_{2}$ into $8$-cycles given by: $ C^{1}=(a_{1}b_{2}a_{3}b_{4}a_{5}b_{6}a_{7}b_{0}),\ C^{2}=(b_{1}a_{2}b_{3}a_{4}b_{5}a_{6}b_{7}a_{0}),\ C^{3}=(a_{1}a_{2}a_{3}a_{4}a_{5}a_{6}a_{7}a_{0})$ and $C^{4}=(b_{1}b_{2}b_{3}b_{4}b_{5}b_{6}b_{7}b_{0}) $. The required $(8; p, q)$-decomposition of $C_{8} \circ \overline{K}_{2}$ is as follows:\\
	$D_{1}$: Consider $C^{1}\cup C^{3}$ and decompose them into $P_{8}$ as follows: $ P^{1}=a_{0}a_{1}b_{2}a_{3}b_{4}a_{5}b_{6}a_{7}b_{0} $ and $ P^{2}=b_{0}a_{1}a_{2}a_{3}a_{4}a_{5}a_{6}a_{7}a_{0} $. \\
	$D_{2}$: Decomposition $C^{2}\cup C^{4}$ into $P_{8}$ given by: $ P^{3}=b_{0}b_{1}a_{2}b_{3}b_{4}b_{5}b_{6}b_{7}a_{0} $ and $ P^{4}=a_{0}b_{1}b_{2}b_{3}a_{4}b_{5}a_{6}b_{7}b_{0} $.\\
	$D_{3}$: Consider $C^{2}\cup C^{3}\cup C^{4}$ and decompose them into $P_{8}$, say $ P^{4}=a_{0}b_{1}b_{2}b_{3}a_{4}\\b_{5}a_{6}b_{7}b_{0},\ P^{5}=b_{0}b_{1}a_{2}a_{1}a_{0}a_{7}a_{6}a_{5}a_{4} $ and $ P^{6}=a_{4}a_{3}a_{2}b_{3}b_{4}b_{5}b_{6}b_{7}a_{0} $.\\
	Here $ (p,q) \in \{(0,4),(2,2),(3,1),(4,0)\}$. For $1 \leq i \leq 3$, combinations of $D_{i}$ provide the required $(8; p, q)$-decomposition.\hfill$\Box$
	
	\begin{lemma}
		For any odd integer $n \geq 3$, $C_{8} \circ \overline{K}_{n}$ admits a $(8; p, q)$-decomposition.
		\label{LEMC_{8}K_{3}C8P8}
	\end{lemma}
	\noindent{\bf Proof.}
	As $n$ is odd, it can be expressed as $3+2t$ for some integer $t \geq 0,\ C_{8} \circ \overline{K}_{n} \cong (C_{8} \circ \overline{K}_{3}) \oplus \left((\frac{n-3}{2})C_{8} \circ \overline{K}_{2}\right)\oplus \left((\frac{n^{2}-2n-3}{4}) P_{2} \times C_{8}\right)$. The graph $ C_{8} \circ \overline{K}_{2}$ has a $(8; p, q)$-decomposition by Lemma \ref{LEMC_{8}K_{2}C8P8}. Let the vertices of $C_{8} \circ \overline{K}_{3}$ be partitioned into $A=\{a_{i}| 0 \leq i \leq 7\},\ B=\{b_{i}| 0 \leq i \leq 7\}$ and $C=\{c_{i}| 0 \leq i \leq 7\}$. The decomposition of $C_{8} \circ \overline{K}_{3}$ into $8$-cycles given by: $ C^{1}=(a_{0}a_{1}a_{2}a_{3}a_{4}a_{5}a_{6}a_{7}),\ C^{2}=(c_{0}a_{1}c_{2}a_{3}c_{4}a_{5}c_{6}a_{7}),\ C^{3}=(b_{0}b_{1}b_{2}b_{3}b_{4}b_{5}b_{6}b_{7}),\ C^{4}=(b_{0}a_{1}b_{2}a_{3}b_{4}a_{5}b_{6}a_{7}),\ C^{5}=(c_{0}c_{1}c_{2}c_{3}c_{4}c_{5}c_{6}c_{7}),\ C^{6}=(c_{0}b_{1}c_{2}b_{3}c_{4}b_{5}c_{6}b_{7}),\ C^{7}=(a_{0}c_{1}a_{2}c_{3}a_{4}c_{5}a_{6}c_{7}),\ C^{8}=(a_{0}b_{1}a_{2}b_{3}a_{4}b_{5}a_{6}b_{7})$ and $C^{9}=(b_{0}c_{1}b_{2}c_{3}b_{4}c_{5}b_{6}c_{7}) $. \\
	  Let the vertices of $P_{2} \times C_{8}$ be partitioned into $D=\{d_{i}| 0 \leq i \leq 7\},\ E=\{e_{i}| 0 \leq i \leq 7\}$ and $F=\{f_{i}| 0 \leq i \leq 7\}$. 	The decomposition of $P_{2} \times C_{8}$ into $8$-cycles given by: $C^{10}=(d_{0}e_{1}d_{2}e_{3}d_{4}e_{5}d_{6}e_{7}),\ C^{11}=(f_{0}e_{1}f_{2}e_{3}f_{4}e_{5}f_{6}e_{7}),\ C^{12}=(e_{0}d_{1}e_{2}d_{3}e_{4}d_{5}e_{6}d_{7})$ and $C^{13}=(e_{0}f_{1}e_{2}f_{3}e_{4}f_{5}e_{6}f_{7})$. The required $(8; p, q)$-decomposition is as follows:\\
	 	$D_{1}$: Since $C^{1}\cup C^{2} \cong C^{3}\cup C^{4} \cong C^{5}\cup C^{6},$ it is enough to consider $P_{8}$-decomposition of $C^{1}\cup C^{2}$. The $8$-paths of $C^{1}\cup C^{2}$ are $ P^{1}=c_{0}a_{1}a_{2}a_{3}a_{4}a_{5}a_{6}a_{7}a_{0} $ and $ P^{2}=a_{0}a_{1}c_{2}a_{3}c_{4}a_{5}c_{6}a_{7}c_{0} $.\\
	 $D_{2}$: Consider $C^{7}\cup C^{8}$ and decompose them into $P_{8}$ as follows: $ P^{3}=b_{7}a_{0}c_{1}a_{2}c_{3}a_{4}c_{5}a_{6}c_{7} $ and $ P^{4}=c_{7}a_{0}b_{1}a_{2}b_{3}a_{4}b_{5}a_{6}b_{7} $. \\
	  $D_{5}$: Consider $C^{7}\cup C^{8}\cup C^{9}$ and decompose them into $P_{8}$ as follows: $ P^{9}=b_{0}c_{7}a_{0}c_{1}b_{2}c_{3}a_{4}c_{5}b_{4},\ P^{10}=b_{4}c_{3}a_{2}b_{1}a_{0}b_{7}a_{6}c_{5}b_{6} $ and $ P^{11}=b_{6}c_{7}a_{6}b_{5}a_{4}b_{3}a_{2}c_{1}b_{0} $. \\
	  Since there are $(\frac{n^{2}-2n-3}{4}) P_{2} \times C_{8}$, here we give $P_{8}$-decomposition of $P_{2} \times C_{8}$ from obtained $C_{8}$-decompositions $\{C^{10},C^{11},C^{12},C^{13}\}$.\\
	$D_{3}$: Consider $C^{10}\cup C^{11}$ and decompose them into $P_{8}$ as follows:$ P^{5}=f_{0}e_{1}d_{2}e_{3}d_{4}e_{5}d_{6}e_{7}d_{0} $ and $ P^{6}=d_{0}e_{7}d_{6}e_{5}d_{4}e_{3}d_{2}e_{1}f_{0} $.\\
	$D_{4}$: Consider $C^{12}\cup C^{13}$ and decompose them into $P_{8}$ as follows:  $ P^{7}=d_{1}e_{0}f_{1}e_{2}f_{3}e_{4}f_{5}e_{6}f_{7} $ and $ P^{8}=f_{7}e_{0}d_{7}e_{6}d_{5}e_{4}d_{3}e_{2}d_{1} $. \\
		Here $(p,q) \in \{(0,13),(2,11),(3,10),(4,9),(5,8),(6,7),(7,6),(8,5),(9,4),\\(10,3),(11,2),(12,1),(13,0)\}$. For $1 \leq i \leq 5$, combinations of $D_{i}$ provide the required $(8; p, q)$-decomposition.\hfill$\Box$
	
		\begin{lemma}
		For any odd integer $n \geq 3$ and $m \equiv 0,1($mod $16)$, $K_{m} \circ \overline{K}_{n}$ admits a $(8; p, q)$-decomposition.
		\label{LEMK_{16}K_{n}C8P8}
	\end{lemma}
	\noindent{\bf Proof.}
	When $m \equiv 1($mod $16)$, $K_{m}$ is decomposed into $C_{8}$ by Theorem \ref{THM5C8P8}. Hence by Lemma \ref{LEMC_{8}K_{3}C8P8}, $C_{8} \circ \overline{K}_{n}$ has a $(8; p, q)$-decomposition.\\
	When $m \equiv 0($mod $16)$, $K_{m}$ is decomposed into subgraphs isomorphic to $K_{16}$ and $K_{8,8}$. By Lemma \ref{LEM3C8P8}, $K_{8,8}$ has a $(8; p, q)$-decomposition. Let the vertices of $K_{16}$ be partitioned into $A=\{a_{i}| 0 \leq i \leq 7\}$ and $B=\{b_{i}| 0 \leq i \leq 7\}$. The $8$-paths of $K_{16}$ is given by: $ P^{1}=a_{0}a_{1}a_{3}a_{6}a_{2}a_{7}a_{5}a_{4}b_{4},\ P^{2}=b_{0}b_{1}b_{3}b_{6}b_{2}b_{7}b_{5}b_{4}a_{3},\ P^{3}=a_{1}a_{2}a_{4}a_{7}a_{3}a_{0}a_{6}a_{5}b_{5},\ P^{4}=b_{1}b_{2}b_{4}b_{7}b_{3}b_{0}b_{6}b_{5}a_{4},\ P^{5}=b_{2}a_{2}a_{3}a_{5}a_{0}a_{4}a_{1}a_{7}a_{6},\ P^{6}=a_{5}b_{2}b_{3}b_{5}b_{0}b_{4}b_{1}b_{7}b_{6},\ P^{7}=a_{2}b_{3}b_{4}b_{6}b_{1}b_{5}b_{2}b_{0}b_{7}$ and $ P^{8}=b_{3}a_{3}a_{4}a_{6}a_{1}a_{5}a_{2}a_{0}a_{7} $. The remaining edges of $K_{16}$ are decomposed into $C_{8}$ as follows: $ C^{1}=(a_{0}b_{1}a_{1}b_{6}a_{6}b_{7}a_{7}b_{0}),\ C^{2}=(a_{0}b_{5}a_{2}b_{4}a_{7}b_{6}a_{3}b_{2}),\ C^{3}=(a_{0}b_{3}a_{1}b_{2}a_{7}b_{1}a_{4}b_{7}),\ C^{4}=(a_{1}b_{4}a_{5}b_{6}a_{4}b_{3}a_{6}b_{0}),\ C^{5}=(a_{2}b_{1}a_{6}b_{5}a_{3}b_{0}a_{5}b_{7}),\ C^{6}=(a_{1}b_{5}a_{7}b_{3}a_{5}b_{1}a_{3}b_{7})$ and $C^{7}=(a_{2}b_{6}a_{0}b_{4}a_{6}b_{2}a_{4}b_{0}) $.  The required $(8; p, q)$-decomposition is as follows:\\
	$D_{1}$: Consider $C^{1}\cup P^{3}$ and decompose them into $P_{8}$ as follows: $P^{1}=b_{5}a_{5}a_{6}b_{6}a_{1}b_{1}a_{0}b_{0}a_{7} $ and $ P^{2}=a_{7}b_{7}a_{6}a_{0}a_{3}a_{7}a_{4}a_{2}a_{1} $.\\
	$D_{2}$: Consider $C^{1}\cup C^{2}$ and decompose them into $P_{8}$ as follows: $P^{3}=b_{2}a_{0}b_{1}a_{1}b_{6}a_{6}b_{7}a_{7}b_{0} $ and $ P^{4}=b_{0}a_{0}b_{5}a_{2}b_{4}a_{7}b_{6}a_{3}b_{2} $.\\
	$D_{3}$: Consider $C^{3}\cup C^{4}$ and decompose them into $P_{8}$ as follows: $P^{5}=a_{6}b_{3}a_{1}b_{2}a_{7}b_{1}a_{4}b_{7}a_{0} $ and $ P^{6}=a_{0}b_{3}a_{4}b_{6}a_{5}b_{4}a_{1}b_{0}a_{6} $. \\
	$D_{4}$: Consider $C^{5}\cup C^{7}$ and decompose them into $P_{8}$ as follows: $P^{7}=b_{6}a_{2}b_{1}a_{6}b_{5}a_{3}b_{0}a_{5}b_{7} $ and $ P^{8}=b_{7}a_{2}b_{0}a_{4}b_{2}a_{6}b_{4}a_{0}b_{6} $.\\
	$D_{5}$: Consider $C^{3}\cup C^{4}\cup C^{6}$ and decompose them into $P_{8}$ as follows: $P^{6}=a_{0}b_{3}a_{4}b_{6}a_{5}b_{4}a_{1}b_{0}a_{6} ,\ P^{9}=a_{0}b_{7}a_{1}b_{5}a_{7}b_{3}a_{5}b_{1}a_{4} $ and $ P^{10}=a_{4}b_{7}a_{3}b_{1}a_{7}b_{2}a_{1}b_{3}a_{6} $.
		Here $ (p,q) \in \{(8,7),(9,6),(10,5),(11,4),(12,3),(13,2),(14,1),(15,0)\}$. For $1 \leq i \leq 5$, combinations of $D_{i}$ provide the required $(8; p, q)$-decomposition.\hfill$\Box$
	
	\begin{theorem}
		$ (K_{m} \circ \overline{K}_{n})(\lambda)$ admits a $(8; p, q)$-decomposition if and only if $8|\frac{\lambda n^{2}m(m-1)}{2}$ and $mn \geq 9$.
	\end{theorem}
	\noindent{\bf Proof.\\}
	The necessary conditions are obvious and we prove the sufficiency as follows:\\
	
	\noindent{\bf Case 1:} $\lambda=1.$ \\
	\noindent{\bf Subcase 1.1:} $n \equiv 0($mod $4)$.\\
	\noindent Here, $K_{m} \circ \overline{K}_{n} \cong K_{m} \circ \overline{K}_{4k}$. When $k=1$, $m$ must be greater than $2$ to have a $(8;p,q)$-decomposition. It is easy to see that $K_{m} \circ \overline{K}_{4}$ can be decomposed into subgraphs isomorphic to $P_{2} \circ \overline{K}_{4}$ and $ \cup_{i=1}^{3} \langle A_{1}, B_{i}\rangle$ mentioned in Lemma \ref{LEM11C8P8}, hence it admits a $(8; p, q)$-decomposition by Lemmas \ref{LEM2C8P8} and \ref{LEM11C8P8}. When $k \geq 2$, let $K_{m} \circ \overline{K}_{4k} \cong \frac{m(m-1)}{2}  K_{4k,4k}$. Hence, $K_{4k,4k}$ admits a $(8;p,q)$-decomposition by Lemma \ref{LEM11C8P8}. \\
	
	\noindent{\bf Subcase 1.2:} $n \equiv 0($mod $2)$ and $m \equiv 0($mod $4)$.\\
	\noindent Here $K_{m} \circ \overline{K}_{n} \cong K_{4k} \circ \overline{K}_{2t}$. By Lemma \ref{LEM00C8P8}, it is enough to consider $K_{4k} \circ \overline{K}_{2}$ which is isomorphic to $K_{8k}-F_{1}$.	When $k$ is even, $  K_{8k}-\textrm{F}_{1} \cong \frac{k}{2} (K_{16}- \textrm{F}_{1}) \oplus \frac{k(k-2)}{2} K_{8,8} $. When $k$ is odd, $ K_{8k}-\textrm{F}_{1} \cong \frac{k-1}{2} (K_{16}- \textrm{F}_{1}) \oplus K_{8}- \textrm{F}_{1}  \oplus \frac{(k-1)^{2}}{2}  K_{8,8} $. We choose one $K_{8,8}$ to be concatenated with $K_{8}- \textrm{F}_{1}$ such that the resultant graph can be expressed in terms of graph considered in Lemma \ref{LEM9C8P8}. Hence it admits a $(8;p,q)$-decomposition. By Lemma \ref{LEM3C8P8} and \ref{LEM8C8P8}, $K_{8,8}$ and $K_{16}-F_{1}$ admits a $(8;p,q)$-decomposition. \\
	
	\noindent{\bf Subcase 1.3:} $n \equiv 0($mod $2)$ and $m \equiv 1($mod $4)$.\\	
	\noindent Clearly $K_{m} \circ \overline{K}_{n}\cong K_{4k+1} \circ \overline{K}_{2t}$. By Lemma \ref{LEM00C8P8}, it is enough to decompose $K_{4k+1} \circ \overline{K}_{2}$. Let $K_{4k+1} \circ \overline{K}_{2} \cong ( k   K_{5}\oplus \frac{k(k-1)}{2}  K_{4,4})\circ \overline{K}_{2} $. Now $ K_{4,4} \circ\overline{K}_{2}$ can be expressed in terms of $P_{2} \circ \overline{K}_{4}$ which admits a $(8;p,q)$-decomposition by Lemma \ref{LEM2C8P8}. As $K_{5} \circ\overline{K}_{2}\cong K_{10}-F_{1}$, it admits a $(8;p,q)$-decomposition by Lemma \ref{LEM12C8P8}. \\
	
	\noindent{\bf Subcase 1.4:} $n \equiv 1($mod $2)$ and $m \equiv 0,1($mod $16)$.\\
	\noindent By Lemma \ref{LEMK_{16}K_{n}C8P8}, $K_{m} \circ \overline{K}_{n}$ admits a $(8;p,q)$-decomposition.\\	
	
	\noindent Now $(8;p,q)$-decomposition of $ K_{m} \circ \overline{K}_{n}(\lambda)$ (where $\lambda$ is any non negative integers) follows for Case 1.\\
	
	\noindent{\bf Case 2:} $\lambda=2$ .\\
	\noindent{\bf Subcase 2.1:} $n \equiv 1($mod $2)$ and $m \equiv  0,1$(mod $8$).\\
 	\noindent By Theorem \ref{THM5C8P8}, $ K_{8k}(2)$ and $ K_{8k+1}(2)$ have a $C_{8}$-decomposition. The proof follows from Lemma \ref{LEMC_{8}K_{3}C8P8}.  \\
	
	\noindent{\bf Subcase 2.2:} $n \equiv  0$(mod $2$).\\
	\noindent By Lemma \ref{LEM2K_{m}{K}_{2}C8P8}, $(K_{m} \circ \overline{K}_{n})(2)$ admits a $(8;p,q)$-decomposition.\\
	
	\noindent{\bf Case 3:} $\lambda=4$ and $n \equiv 1($mod $2)$ and $m \equiv  0,1$(mod $4$).\\
	\noindent By Theorem \ref{THM5C8P8}, $ K_{4k}(4)$ and $K_{4k+1}(4)$ have a $C_{8}$-decomposition. The proof follows from Lemma \ref{LEMC_{8}K_{3}C8P8}.  \\
	
	\noindent This completes the proof of the theorem.\hfill$\Box$

\end{document}